\documentclass[framed]{customarticle}

\title{From Worst-Case to Smoothed Complexity \\ in Continuous Optimization}
\author{Daniel Cortild \orcidCortild, 
Coralia Cartis \orcidCartis \\ Mathematical Institute, University of Oxford, United Kingdom}
\date{Last compiled: \today}

\usepackage{enumitem}

\begin{document}

\maketitle

\begin{abstract}
\noindent
Worst-case complexity bounds often fail to reflect the practical performance of optimization algorithms, while average-case complexity heavily relies on the distribution of the input. In this work, we develop a smoothed-analysis framework for classical continuous optimization algorithms by studying their behaviour on random perturbations of arbitrary functions. Our main contribution is to apply methods suited for non-degenerate functions to perturbations of degenerate functions, thereby recovering the faster best-known rates for the non-degenerate class. We derive smoothed complexity guarantees for several standard methods, including (Accelerated) Gradient Descent, Conjugate Gradient, and Newton-type methods, under various settings including nonconvex ones, and show that these bounds can substantially improve upon their worst-case guarantees, bringing them closer to the practical behaviour of the respective methods.

\end{abstract}


\section{Introduction}





Many theoretical results quantifying the behaviour  of algorithms for optimization problems are framed as \textit{oracle complexity bounds}, also frequently referred to as \textit{evaluation complexity bounds}. This metric quantifies the number of \textit{oracle calls}, denoted by $\mathbb K^\varepsilon(f)$, required by a given algorithm to converge to an $\varepsilon$-optimal point for a problem characterized by a mapping $f$. An oracle call typically refers to a black-box which, given a point, returns pointwise information about the function, but no global information. Such an oracle could, for instance, return function values, gradient evaluations, Hessian-vector products, or combinations thereof. The notion of $\varepsilon$-optimality may be defined in several different ways, such as absolute or relative distance to solution, function value accuracy, gradient norm, etc. Even with the oracle and notion of optimality settled, the resulting complexity measure can be measured in different ways, of which we expose three different ones below.

\paragraph{Worst-Case Complexity.}
Historically, the vast majority of the optimization literature focuses heavily on \textit{worst-case} oracle complexity guarantees. The goal of such an analysis is to provide a uniform upper bound on the oracle complexity $\mathbb K^\varepsilon(f)$ that remains valid uniformly for all instances $f\in \mathcal F$ within a specified problem class $\mathcal F$ \citep{nemirovskii_problem_1983, vavasis_nonlinear_1991, nesterov_lectures_2018, cartis_evaluation_2022, wright_optimization_2022}. 

Even when these upper bounds are tight, meaning there exists a problem within $\mathcal F$ that attains the upper bound, such instances are often pathological in some way and characterizing the practical performance of the algorithm solely with this worst-case scenario is often overly pessimistic. Such an analysis frequently fails to explain why algorithms perform well on real-world problems that lack this worst-case geometry. 


\paragraph{Average-Case Complexity.}
To address this problem, alternative complexity measures have been proposed. The first is the one of \textit{average-case} complexity, which computes the expectation of the complexity $\ExpD{f\sim \mathcal D}{\mathbb K^\varepsilon(f)}$, where $\mathcal D$ is a probability distribution over $\mathcal F$. By integrating over an entire ensemble of instances, the framework captures the ``typical'' performance that we can expect to encounter in practice. Presenting an average-case complexity result next to a worst-case complexity result allows to gain an image of the pessimistic and the typical behaviour of the algorithm. Famous examples show that quicksort, despite having worst-case complexity $\mathcal O(n^2)$, has average-case complexity $\mathcal O(n\log n)$ \citep{hoare_quicksort_1962}, that quickselect has worst-case complexity $\mathcal O(n^2)$ yet average-case complexity $\mathcal O(n)$ \citep{hoare_algorithm_1961}, and that the simplex method, which can exhibit exponential-time behaviour, runs in polynomial time in average \citep{borgwardt_distributionindependent_1982}.

Optimization algorithms have recently been studied in an average-case perspective. Due to the complexity of the question \textit{``What is a random function?''}, most of the works have tackled the study on random quadratic functions. First-order methods have been analyzed in \cite{scieur_universal_2020,cunha_only_2022,pedregosa_averagecase_2020,paquette_halting_2023}, with \cite{benning_gradient_2024} considering random generalized linear models. Another line of work has considered the study of the average-case behaviour of the Conjugate Gradient Method on various distributions of random quadratics, namely \cite{deift_conjugate_2021,ding_conjugate_2021,paquette_universality_2023}. Slightly tangential, but mention-worthy in our opinion, is the line of work in \cite{paquette_sgd_2021,paquette_dynamics_2021,lee_trajectory_2022,marshall_clip_2025}, who study the high-dimensional limiting dynamics of stochastic methods over a random class of quadratic functions. 


However, average-case analyses are hard to carry out in general due to analytical difficulties of probabilistic iterative processes, and to interpret in practice due to their strong dependency on the input distribution, which may poorly reflect instances observed in practice.



\paragraph{Smoothed Complexity.} \textit{Smoothed analysis} aims to bridge the gap between worst-case and typically-observed complexity by studying the expected performance of an algorithm under a small localized perturbation of a worst-case instance. It was first introduced in \cite{spielman_smoothed_2004}, where the authors provided evidence for the simplex algorithm with standard pivoting rules taking polynomial time to run, whereas the worst-case complexity was exponential. This showed that the worst-case instances were rare instances and could be avoided by a small perturbation. A more recent survey on smoothed analysis for linear programming includes \cite{dadush_friendly_2020} and the references therein.

Beyond linear programming, smoothed analysis has been successfully applied across a wide range of algorithmic settings. This includes local-search heuristics for combinatorial problems \cite{manthey_smoothed_2020}, size of Pareto sets in multiobjective optimization \cite{roglin_smoothed_2020}, analysis of the condition number of structured real polynomial systems \cite{ergur_probabilistic_2019}, game-theoretic optimization \cite{anagnostides_convergence_2024}, and online metric problems \cite{coester_smoothed_2025}, to name a few. 


In the same spirit, it was proposed in \cite{menon_smoothed_2016} to study the smoothed complexity of optimization algorithms, namely of the Conjugate Gradient Method. They study the complexity of a perturbed quadratic function, where the perturbation is itself a quadratic function, and show the expected complexity is much better than the worst-case complexity over the class of quadratic functions. 

\subsection{Contribution and Motivation}

In this paper, given a function class $\mathcal F$, we consider the smoothed complexity measure 
\[
    \sup_{f\in \mathcal F}\ExpD{W}{\mathbb K^{\varepsilon}(f_{W})}
\]
of the function class $\mathcal F$, where, given $f\colon \R^n\to \R\in \mathcal F$, $f_{W}$ denotes a random perturbation of the function $f$. This notion will be formalized in Section \ref{sec:Results}.

The smoothed complexity count $\sup_{f\in \mathcal F}\ExpD{W}{\K^\varepsilon(f_{W})}$ interpolates between the worst-case and average-case analysis: an adversary first chooses a function $f\in \mathcal F$, and only the perturbation is randomized. This is the optimization analogue of the smoothed analysis in \cite{spielman_smoothed_2004}, inspired by \cite{menon_smoothed_2016}.

Within this framework, we study Gradient Descent (Section \ref{sec:GD1L}), Nesterov Accelerated Gradient Descent (Section \ref{sec:NAG}), the Conjugate Gradient method (Section \ref{sec:CG}), the Cubic Regularized Newton Method (Section \ref{sec:CNM}), a Hybrid Gradient--Cubic Newton method (Section \ref{sec:Hyb}), and the Gauss--Newton method in the convex and nonconvex settings (Sections \ref{sec:GNConvex}-\ref{sec:GNConvexArmijo} and \ref{sec:GN_PL}, respectively). This is by no means an exhaustive list of the algorithms that can be covered in our setting, but they provide, in our opinion, a representative overview of main optimization methods and of our analytical framework. 

In all the cases we consider, we show that degenerate problem instances benefit, under any positive perturbation, from the same complexity guarantees as the corresponding non-degenerate class. This therefore explains why algorithms perform well in practice, even in degenerate settings. Specifically we observe the following:
\begin{itemize}
    \item For the standard first-order methods of Gradient Descent (Sections \ref{sec:GD1L}-\ref{sec:GD1L_PL}), Nesterov Accelerated Gradient Descent (Section \ref{sec:NAG}), and Conjugate Gradient (Section \ref{sec:CG}), it is known that convex functions enjoy a non-asymptotic complexity guarantee polynomial in $1/\varepsilon$, whereas strongly convex functions improve this to $\log(1/\varepsilon)$. The latter comes at the cost of a dependency on the condition number $\kappa$, which is infinite in the merely convex case, and can be arbitrarily large for function with a very small strong convexity parameter. The same holds true in a nonconvex setting, under gradient domination inequalities. We show that the random perturbation of the input does not make the condition number arbitrarily large in expectation, and that the smoothed complexity guarantee for convex functions enjoys a favourable dependency on $\varepsilon$ without a bad condition number.
    \item For second-order methods, such as Cubic Regularized Newton (Sections \ref{sec:CNM}-\ref{sec:CNM_PL}) and a hybrid method between Gradient Descent and Cubic Regularized Newton (Section \ref{sec:Hyb}), classical global rate bounds show a constant number of steps followed by local quadratic convergence on strongly convex functions. We show that convex functions, under an arbitrary positive perturbation, enjoy the same type of non-asymptotic convergence bound, and quantify the size of the quadratic convergence region and the dependence on the problem condition number. We also study a nonconvex setting under gradient domination.
    \item Finally, for Gauss-Newton methods (Section \ref{sec:GNConvex}-\ref{sec:GN_PL}), we show that perturbations allow the method to retain linear convergence in the convex case or radially monotone case. We also study an adaptive Armijo line-search variant in the convex case.
\end{itemize}

\paragraph{Motivation.} As a secondary motivation to computing the quantity in Equation \eqref{eq:to_be_computed}, we mention the recent work on differential privacy through additive random quadratics in \cite{cortild_quadratic_2026}. In this work, it is shown that an additive quadratic perturbation inherently enforces privacy when solving the problem, without loosing too much on the accuracy side. No complexity results were provided, as the study of such complexity notion was out of  the scope of the cited paper. The current paper complements that story, by providing the necessary complexity computations.

\paragraph{Notation.} We denote by, for any base $b> 1$, $\log_{b,+}(x)=\log_b(x)$ if $x\ge 1$ and $0$ otherwise. We note that $\log_{b, +}(x)\le \log_b(1+x)$ for all $x\ge 0$.

\section{Random Matrix Theory Preliminaries}

Consider a random real $n\times m$ Gaussian matrix $G$ (with entries independently following the standard normal distribution), and let $W=\frac1m GG^T$ be a real $n\times n$ Wishart matrix with parameter $m$. We place ourselves in the setting of \textit{proportional growth}, in which the parameters $m$ and $n$ grow proportionally with $n/m\to d\in (0, 1)$. We assume this limit is taking place from the left, so that $n/m\le d$ holds for all $(n, m)$. We note that this corresponds to the \textit{overparametrized model}, such that $W$ is almost surely nonsingular. For a fixed $n$, as $d\to 0$, we must have $m\to \infty$, corresponding to the infinitely overparametrized setting.

The results in this paper reduce to the study of the extreme eigenvalues of $W$, denoted by $\lambda_{\min}(W)$ and $\lambda_{\max}(W)$, Equivalently, we study the extreme singular values of $G$, denoted by $s_{\min}(G)$ and $s_{\max}(G)$, linked to the extreme eigenvalues of $W$ through the relations
\[
    \lambda_{\min}(W)=\frac{s_{\min}(G)^2}{m}\quad\text{and}\quad \lambda_{\max}(W)=\frac{s_{\max}(G)^2}{m}.
\]
Wishart matrices and its singular values have been thoroughly studied over the past decades. To mention a few, we refer to \cite{bai_limit_1993, davidson_local_2001, johnstone_distribution_2001, dumitriu_matrix_2002, edelman_tails_2005, rudelson_smallest_2009, rudelson_nonasymptotic_2010, chiani_distribution_2014,huang_remark_2020}. For a more recent survey and introduction to the matter we refer to \cite{vershynin_highdimensional_2018} and the references therein.

We now introduce the main moment bounds used throughout. They can be obtained from, for instance, \cite[Lemma 2.2, Theorem 2.2]{sarkar_moment_2026}. However, the resulting bounds are quite loose, so we reprove them. The proofs and explicit constants are given in Appendix \ref{app:constants}.

\begin{theorem}[Moment Bounds of Singular Values]\label{thm:moments}\label{coro:mixed_moments}
    For every $q>0$, there exists constants $C^{\max}_{q, d}$ and $C^{\min}_{d, q}$ so that, for $m, n$ sufficiently large,
    \[
        \Exp{s_{\max}(G)^qm^{-q/2}}\le C^{\max}_{d, q}\quad \text{and}\quad \Exp{s_{\min}(G)^{-q}m^{q/2}}\le C^{\min}_{d, q}.
    \]
\end{theorem}

\section{Smoothed Analysis Results}\label{sec:Results}

In this paper, we consider the smoothed complexity measure over a function class $\mathcal F$ given by
\begin{equation}\label{eq:to_be_computed}
    \sup_{f\in \mathcal F}\ExpD{W}{\mathbb K^\varepsilon(f_W)}, \quad \text{where $f_W(x)=f(x)+\frac{\sigma^2}{2}(x-\tilde x)^TW(x-\tilde x)$},   
\end{equation}
where $f\colon \R^n\to \R\in \mathcal F$ is a function belonging to the function class $\mathcal F$,
and $f_W$ is a perturbation of the function $f$ characterized by the scaling parameter $\sigma$, a given center $\tilde x$ and the random Wishart matrix $W=\frac{1}{m}GG^T$.\footnote{Equivalently, we could define $f_{W}(x)=f(x)+\frac{1}{2}(x-\tilde x)^TW(x-\tilde x)$, where $W=\frac{\sigma^2}{m}GG^T$.} The scaling parameter $\sigma$ plays the role of the variance of the random matrix $W$. We note that $\sigma, x_0$ and $\tilde x$ are fixed before drawing $W$, and are fixed accross the supremum.

All of our guarantees will concern optimization of the perturbed objective $f_W$. When $f$ has a minimizer $x_*$, guarantees on $f_W$ are linked to guarantees on $f$ through 
\[
    f(x)-\min f\le \delta+\frac{\sigma^2}{2}(x_*-\tilde x)^TW(x_*-\tilde x),
\]
where $f_W(x)-\min f_W\le \delta$. Thus, the improved regularity comes with a perturbation bias.



For convex $L$-smooth $f$, it holds that $f_W$ is $\mu_\sigma$-strongly convex and $L_\sigma$-smooth, for
\begin{equation}\label{eq:mu_L}
    \mu_\sigma=\sigma^2\lambda_{\min}(W)=\sigma^2\frac{s_{\min}(G)^2}{m}\quad \text{and}\quad L_\sigma= L+\sigma^2\lambda_{\max}(W)=L+\sigma^2\frac{s_{\max}(G)^2}{m}.
\end{equation}
For $D^2=\|x_0-\tilde x\|^2$ and $\Delta=f(x_0)-\inf f$,
\begin{equation}\label{eq:Del0}
    f_W(x_0)-\inf f_W\le \Delta_\sigma=\Delta+\frac{\sigma^2 D^2\lambda_{\max}(W)}2=\Delta+\frac{\sigma^2 D^2s_{\max}(G)^2}{2m}.
\end{equation}

\textbf{Numerical Illustrations.} For each method, we compare
\begin{itemize}
    \item the theoretical smoothed complexity bound;
    \item the empirical expectation of the corresponding worst-case bound; and
    \item the empirical iteration count of the algorithm when applied to $f_W$.
\end{itemize}
All empirical quantities are computed as an average over $100$ samples of the perturbation. The experiments were implemented in Python 3.13 and run on a MacBook Pro with an Apple M5 chip.\footnote{Code available on \url{https://github.com/DanielCortild/Smoothed-Analysis}.} The function used in each experiment is specified in the corresponding section.

\subsection{Gradient Descent}\label{sec:GD1L}

We begin with the \textit{Gradient Descent} method \citep{cauchy_methode_1847} on the class $\mathcal F_\mu^L$ of $\mu$-strongly convex and $L$-smooth functions, with the convention that $\mu=0$ corresponds to the convex bounded-below setting $\mathcal F_0^L\supset \mathcal F_\mu^L$. For $h\in \mathcal F_{\mu}^L$, we consider, for $x_0\in \R^n$ and $k\ge 0$,
\begin{equation}\tag{GD}\label{alg:GD1L}
    x_{k+1}=x_k-\frac1L\nabla h(x_k).
\end{equation}
We measure complexity in terms of relative function value accuracy, namely, for $\varepsilon\in (0, 1)$,
\[
    \mathbb K^\varepsilon_{\text{GD}}(h)=\min_{k\ge 0}\left\{k\colon \frac{h(x_k)-\inf h}{h(x_0)-\inf h}\le \varepsilon\right\},
\]
where $\mathbb K^\varepsilon_{\text{GD}}(h)=0$ if $h(x_0)=\inf h$.

We first recall the classical worst-case complexity (see, for instance, \cite[Equation (3.15)]{wright_optimization_2022}), which will serve as a benchmark for our smoothed complexity guarantees.

\begin{theorem}[Worst-Case Complexity for \eqref{alg:GD1L} - {\cite{wright_optimization_2022}}]\label{thm:WC_GD1L}
    Let $h\in\mathcal F_\mu^L$ with $0<\mu\le L<+\infty$, and let $(x_k)$ be the iterates of \eqref{alg:GD1L}. Then, for $\kappa=L/\mu$,
    \[
        \frac{h(x_k)-\inf h}{h(x_0)-\inf h}\le\left(1-\frac1\kappa\right)^k.
    \]
    Thus the relative accuracy complexity is bounded by
    \[
        \sup_{h\in\mathcal F_\mu^L}\K^\varepsilon_{\text{GD}}(h)\le\left\lceil\kappa\log\frac1\varepsilon\right\rceil.
    \]
\end{theorem}

Though the above theorem is only valid for $\mu>0$, it is known that Algorithm \eqref{alg:GD1L} converges for $\mu=0$, albeit at a sublinear rate. 
The below represents our first \textit{smoothed complexity} bound.


\begin{theorem}[Smoothed Complexity for \eqref{alg:GD1L}]\label{thm:SA_GD1L}
    Let $0<\sigma,L<+\infty$ and $\varepsilon\in(0,1)$. We denote by $f_W$ a perturbation of $f$, parametrized by $\sigma$, as defined in Equation \eqref{eq:to_be_computed} and apply Algorithm \eqref{alg:GD1L} to $f_W$. 
    Then, as $n\to\infty$, the smoothed complexity of Algorithm \eqref{alg:GD1L} over the function class $\mathcal F_0^L$ is bounded by
    \[
        \sup_{f\in\mathcal F_0^L}\ExpD{W}{\K^\varepsilon_{\text{GD}}(f_W)}
        \le1+ \mathcal O_d\left((1+\sigma^{-2}L)\log\frac1\varepsilon\right).
    \]
\end{theorem}
\begin{proof}
    Fix any $f\in\mathcal F_0^L$.
    By Equation \eqref{eq:mu_L}, Cauchy--Schwarz and Theorem \ref{thm:moments},
    \[
        \ExpD{W}{\kappa_\sigma}
        =L\sigma^{-2}\ExpD{W}{\frac{m}{s_{\min}(G)^2}}
        +\ExpD{W}{\frac{s_{\max}(G)^2}{s_{\min}(G)^2}}
        \le\mathcal O_d(1+\sigma^{-2}L).
    \]
    The result follows from Theorem \ref{thm:WC_GD1L} applied to $h=f_W$, and where the ensuing step-size in Algorithm \eqref{alg:GD1L} is $1/L_\sigma$, where $L_\sigma$ is defined in \eqref{eq:mu_L}.\qedhere
\end{proof}

\textbf{Discussion.} The smoothed complexity bound depends only logarithmically on $\varepsilon$, even though we are considering the class of merely convex functions, for which we expect an inverse dependency in $\varepsilon$. This is due to the perturbation, which renders the objective strongly convex.
The strong convexity parameter may however be arbitrarily small, in which case the condition number is arbitrarily large, and the linear bound is worse than the sublinear bound. Given the nature of the perturbation matrix, this occurs with sufficiently small probability, so that, in expectation, we do not observe this. In fact, the expected strong convexity parameter is proportional to the scale of the perturbation, $\sigma^2$, as can be seen  in the bound, where $\sigma^{-2}L$ acts as a condition number. 

\textbf{Numerical Illustration.} We illustrate Theorem \ref{thm:SA_GD1L} on a weighted Huber objective, based on the Huber function used to generate bad instances for smooth convex first-order optimization \citep{taylor_smooth_2017},
\[
    f(x)=\sum_{i=1}^n 2^{10(i-n/2)/n}\phi(x_i),
    \qquad
    \phi(t)=\begin{cases}
        \frac12t^2, & |t|\le 1,\\
        |t|-\frac12, & |t|>1.
    \end{cases}
\]
The weights give $L=2^5$ and a ratio of largest to smallest curvature of $2^{10(n-1)/n}$ in the quadratic region.
The results are shown in Figure \ref{fig:GD1L}, including the iteration count of GD applied directly to the unperturbed objective with step-size $1/L$, shown as a dotted line without markers. This additional curve is independent of $\sigma$ and $d$. We observe that our bounds are very close to the empirical smoothed complexity of the worst-case bound, but that both are overly pessimistic compared to the empirical smoothed complexity of the run algorithm. This highlights the well-known fact that first-order methods can be overly pessimistic in their worst-case behaviour. 

\begin{figure}[H]
    \centering
    \def\SAMethod{GD1L}
    \input{GeneratePlot}
    \caption{Smoothed complexity of Algorithm \eqref{alg:GD1L} under Wishart perturbations.}
    \label{fig:GD1L}
\end{figure}

We additionally consider a degenerate ill-conditioned quadratic, given by 
\[
    f(x)=\frac12\sum_{i=1}^{\lfloor n/2\rfloor}a_i x_i^2,
    \qquad a_i=2^{10(i-n/2)/n}.
\]
It has a large spread in its non-zero eigenvalues, and a large number of zero eigenvalues. The function is therefore not strongly convex. The results are shown in Figure \ref{fig:GD1LQuadratic}. We see that our results correctly predict the tendency in $\sigma$, whereas the dependency in $d$ is more pessimistic.

\begin{figure}[H]
    \centering
    \def\SAMethod{GD1LQuadratic}
    \input{GeneratePlot}
    \caption{Smoothed complexity of Algorithm \eqref{alg:GD1L} under Wishart perturbations.}
    \label{fig:GD1LQuadratic}
\end{figure}

\subsection{Gradient Descent under Radial Monotonicity}\label{sec:GD1L_PL}

We now consider Algorithm \eqref{alg:GD1L} beyond the convex setting. Specifically, we say that a function $h\in \mathcal F_{-\infty}^{L}$, where $\mathcal F_{-\infty}^{L}$ denotes the class of nonconvex $L$-smooth functions, satisfies \textit{gradient domination}, also known as the \textit{Polyak-Lojasiewicz inequality}, with parameter $\eta>0$ if it is bounded below and
\[
    \|\nabla h(x)\|^2\ge 2\eta (h(x)-\inf h)\quad \text{for all } x\in \R^n.
\]
Denote by $\mathcal F_{\text{PL}, \eta}^{L}$ the subclass of $\mathcal F_{-\infty}^{L}$ that satisfies gradient domination with parameter $\eta>0$. We note that gradient domination implies that all local minimizers are global, and therefore we can measure complexity in terms of relative function value accuracy $\mathbb K^\varepsilon_{\text{GD}}(h)$ as done in Section \ref{sec:GD1L}.

Gradient domination does not require convexity, but still ensures linear convergence of Gradient Descent. We first give the deterministic worst-case complexity under gradient domination, whose proof may be found in \cite[Theorem 1]{karimi_linear_2020}.

\begin{theorem}[{Worst-Case Complexity of \eqref{alg:GD1L} under Gradient Domination - \cite{karimi_linear_2020}}]\label{thm:WC_GD1L_GradDom}
    Let $h\in \mathcal F_{\text{PL}, \eta}^{L}$ with $0<\eta\le L<+\infty$, and let $(x_k)$ be the iterates of \eqref{alg:GD1L}. Then,
    \[
        \frac{h(x_k)-\inf h}{h(x_0)-\inf h}\le \left(1-\frac{\eta}{L}\right)^k.
    \]
    Thus, for $\varepsilon\in (0, 1)$,. the relative function value accuracy complexity is bounded by 
    \[
        \sup_{h\in \mathcal F_{\text{PL}, \eta}^{L}}\mathbb K^\varepsilon_{\text{GD}}(h)\le \left\lceil\frac{L}{\eta}\log\frac{1}{\varepsilon}\right\rceil.
    \]
\end{theorem}
We now consider the class of \textit{radially monotone} functions, denoted by $\mathcal F_{\text{rad}}^{L}(\tilde x)$, within $\mathcal F_{-\infty}^{L}$, namely functions for which
\[
    f(\tilde x)=\inf f\quad \text{and}\quad (x-\tilde x)^T\nabla f(x)\ge 0\quad \text{for every $x\in \R^n$}.
\]
These functions are nondecreasing along rays to the specified center, which is a minimizer, but could still have negative curvature, and thus need not be convex. 

Fixing the center of the perturbation to be the same $\tilde x$, we obtain that
\[
    (x-\tilde x)^T\nabla f_W(x)=(x-\tilde x)^T\nabla f(x)+(x-\tilde x)^T\sigma^2W(x-\tilde x)\ge \sigma^2\lambda_{\min}(W)\|x-\tilde x\|^2,
\]
which implies that 
\[
    \|\nabla f_W(x)\|\ge \sigma^2\lambda_{\min}(W)\|x-\tilde x\|.
\]
As $f_W$ is $L_\sigma$-smooth for $L_\sigma=L+\sigma^2\lambda_{\max}(W)$, and $\tilde x$ is a minimizer of $f_W$, we get 
\[
    f_W(x)-f_W(\tilde x)\le \frac{L_\sigma}{2}\|x-\tilde x\|^2,
\]
which, combined with the previous point, gives
\[
    \|\nabla f_W(x)\|^2\ge \sigma^4\lambda_{\min}(W)^2\|x-\tilde x\|^2\ge \frac{2\sigma^4\lambda_{\min}(W)^2}{L_\sigma}(f_W(x)-\inf f_W).
\]
Specifically, $f_W\in \mathcal F_{\text{PL}, \eta_\sigma}^{L_\sigma}$, for
\[
    L_\sigma=L+\sigma^2\lambda_{\max}(W)\quad\text{and} \quad \eta_\sigma=\frac{\sigma^4\lambda_{\min}(W)^2}{L_\sigma}.
\]
We now state the smoothed complexity guarantee.
\begin{theorem}[Smoothed Complexity of \eqref{alg:GD1L} under Radial Monotonicity]\label{thm:SA_GD1L_PL}
    Let $0<\sigma, L <+\infty$, $\varepsilon\in (0, 1)$, and $\tilde x\in \R^n$. We denote by $f_W$ a perturbation of $f$, parametrized by $\sigma$, as defined in Equation \eqref{eq:to_be_computed}. Then, as $n\to \infty$, the smoothed complexity of Algorithm \eqref{alg:GD1L} over the function class $\mathcal F_{\text{rad}}^{L}(\tilde x)$ is bounded by 
    \[
        \sup_{f\in \mathcal F_{\text{rad}}^{L}(\tilde x)}\ExpD{W}{\mathbb K^\varepsilon_{\text{GD}}(f_W)}\le 1+ \mathcal O_d\left((1+\sigma^{-2}L)^2\log\frac{1}{\varepsilon}\right).
    \]
\end{theorem}
\begin{proof}
    Fix any $f\in \mathcal F_{\text{rad}}^{L}(\tilde x)$, so that $f_W\in \mathcal F_{\text{PL}, \eta_\sigma}^{L_\sigma}$. Theorem \ref{thm:WC_GD1L_GradDom} and the preceding bounds give 
    \[
        \mathbb K^\varepsilon_{\text{GD}}(f_W)\le1+ \frac{(L+\sigma^2\lambda_{\max}(W))^2}{\sigma^4\lambda_{\min}(W)^2}\log\frac{1}{\varepsilon}.
    \]
    Expanding the numerator, applying Cauchy-Schwarz to each mixed moment and using Theorem \ref{thm:moments} proves the result.
\end{proof}

\textbf{Discussion.} As in Section \ref{sec:GD1L}, the smoothed complexity bound depends only logarithmically on $\varepsilon$, even though the original objective did not satisfy gradient domination. In contrast to Section \ref{sec:GD1L}, the dependency in $1+\sigma^{-2}L$ is squared. This is explained through the equality $\frac{L_\sigma}{\eta_\sigma}=\frac{L_\sigma^2}{\mu_\sigma^2}$, where $\mu_\sigma=\sigma^2\lambda_{\min}(W)$. This additional factor thus arises naturally in the passage from radial monotonicity to gradient domination. 

\textbf{Numerical Illustration.} To illustrate Theorem \ref{thm:SA_GD1L_PL}, we consider the objective
\[
    f(x)=\sum_{i=1}^n\left(1-e^{-x_i^2/2}\right), \quad \tilde x=0.
\]
This objective satisfies $(x-\tilde x)^T\nabla f(x)=\sum_i x_i^2e^{-x_i^2/2}\ge 0$ and $f(\tilde x)=\min f$, so that $f$ is radially monotone. However, $f$ is nonconvex, as $\partial_{ii}f(x)=(1-x_i^2)e^{-x_i^2/2}$, and does not satisfy a global gradient domination bound, as, for instance, along $x=te_1$, the function tends to a finite value as $t\to \infty$. The objective is however globally smooth with $L=1$. 

\begin{figure}[H]
    \centering
    \def\SAMethod{GDRadial}
    \input{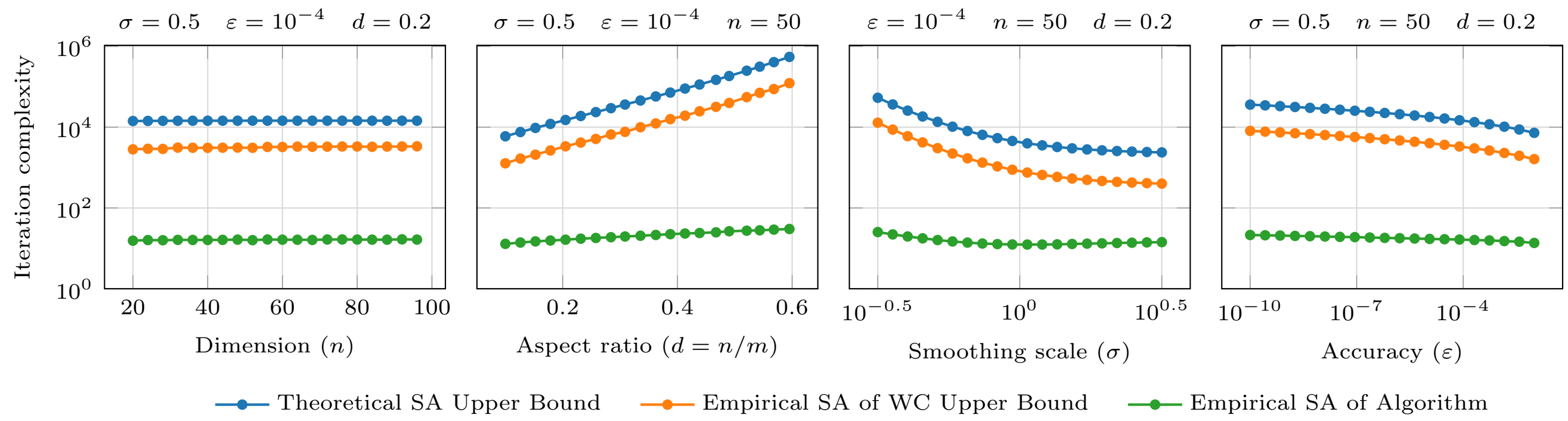}
    \caption{Smoothed complexity of Algorithm \eqref{alg:GD1L} under Wishart perturbations.}
    \label{fig:GD1L_PL}
\end{figure}

The results are shown in Figure \ref{fig:GD1L_PL}. We note that the growth in $d$ and in $\sigma$ is pessimistic in our bounds, but the growth occurs in the correct trend. However, these results show the applicability of the PL guarantee in our framework. 

\subsection{Nesterov Accelerated Gradient Descent}\label{sec:NAG}

We next consider \textit{Nesterov Accelerated Gradient Descent} \citep{nesterov_method_1983} on the same class $\mathcal F_\mu^L$. For $h\in\mathcal F_\mu^L$, initialize $y_0=x_0$ and consider, for $k\ge 0$,
\begin{equation}\tag{NAG}\label{alg:NAG}
    x_{k+1}=y_k-\frac1L\nabla h(y_k),
    \qquad
    y_{k+1}=x_{k+1}+\frac{\sqrt L-\sqrt\mu}{\sqrt L+\sqrt\mu}(x_{k+1}-x_k).
\end{equation}
We again measure complexity in terms of relative function value accuracy, namely, for $\varepsilon\in (0, 1)$,
\[
    \mathbb K^\varepsilon_{\text{NAG}}(h)=\min_{k\ge 0}\left\{k\colon \frac{h(x_k)-\inf h}{h(x_0)-\inf h}\le \varepsilon\right\},
\]
where $\mathbb K^\varepsilon_{\text{NAG}}(h)=0$ if $h(x_0)=\inf h$.

We first recall the classical accelerated worst-case complexity; see, for instance, \cite[Theorem 2.2.1]{nesterov_lectures_2018}. The bound is accelerated in the sense that it depends on $\sqrt{\kappa}$ rather than $\kappa$, as in Section \ref{sec:GD1L}.


\begin{theorem}[Worst-Case Complexity for \eqref{alg:NAG} - {\cite{nesterov_lectures_2018}}]\label{thm:WC_NAG}
    Let $h\in\mathcal F_\mu^L$ with $0<\mu\le L<+\infty$, and let $(x_k)$ be the iterates of \eqref{alg:NAG}. Then, for $\kappa=L/\mu$,
    \[
        \frac{h(x_k)-\inf h}{h(x_0)-\inf h}
        \le 2\left(1-\frac1{\sqrt\kappa}\right)^k.
    \]
    Thus the relative accuracy complexity is bounded by
    \[
        \sup_{h\in \mathcal F_\mu^L}\K^\varepsilon_{\text{NAG}}(h)
        \le \left\lceil\sqrt\kappa\log\frac2\varepsilon\right\rceil.
    \]
\end{theorem}

We provide our smoothed complexity guarantee for Algorithm \eqref{alg:NAG}, which also shows the acceleration behaviour.


\begin{theorem}[Smoothed Complexity for \eqref{alg:NAG}]\label{thm:SA_NAG}
    Let $0<\sigma, L<+\infty$ and $\varepsilon\in(0,1)$. We denote by $f_W$ a perturbation of $f$, parametrized by $\sigma$, as defined in Equation \eqref{eq:to_be_computed}. Then, as $n\to\infty$, the smoothed complexity of Algorithm \eqref{alg:NAG} over the function class $\mathcal F_0^L$ is bounded by
    \[
        \sup_{f\in\mathcal F_0^L}\ExpD{W}{\K^\varepsilon_{\text{NAG}}(f_W)}
        \le 1+\mathcal O_d\left(\left(1+\sigma^{-1}\sqrt L\right)\log\frac2\varepsilon\right).
    \]
\end{theorem}
\begin{proof}
    Fix any $f\in\mathcal F_0^L$. By Equation \eqref{eq:mu_L}, subadditivity of the square root, Cauchy--Schwarz and Theorem \ref{thm:moments},
    \[
        \ExpD{W}{\sqrt{\kappa_\sigma}}
        \le \sigma^{-1}\sqrt L\ExpD{W}{\frac{\sqrt m}{s_{\min}(G)}}
        +\ExpD{W}{\frac{s_{\max}(G)}{s_{\min}(G)}}
        \le \mathcal O_d(1+\sigma^{-1}\sqrt L).
    \]
    The result follows from Theorem \ref{thm:WC_NAG}. \qedhere
\end{proof}

\textbf{Discussion.} As in Section \ref{sec:GD1L}, the smoothed complexity bound depends logarithmically on $\varepsilon$, due to the perturbation inducing strong convexity. Moreover, $\sigma^{-2}L$ acts as a condition number, so that the smoothed complexity bound depends on the square root of the condition number.

\begin{figure}[H]
    \centering
    \def\SAMethod{NAG}
    \input{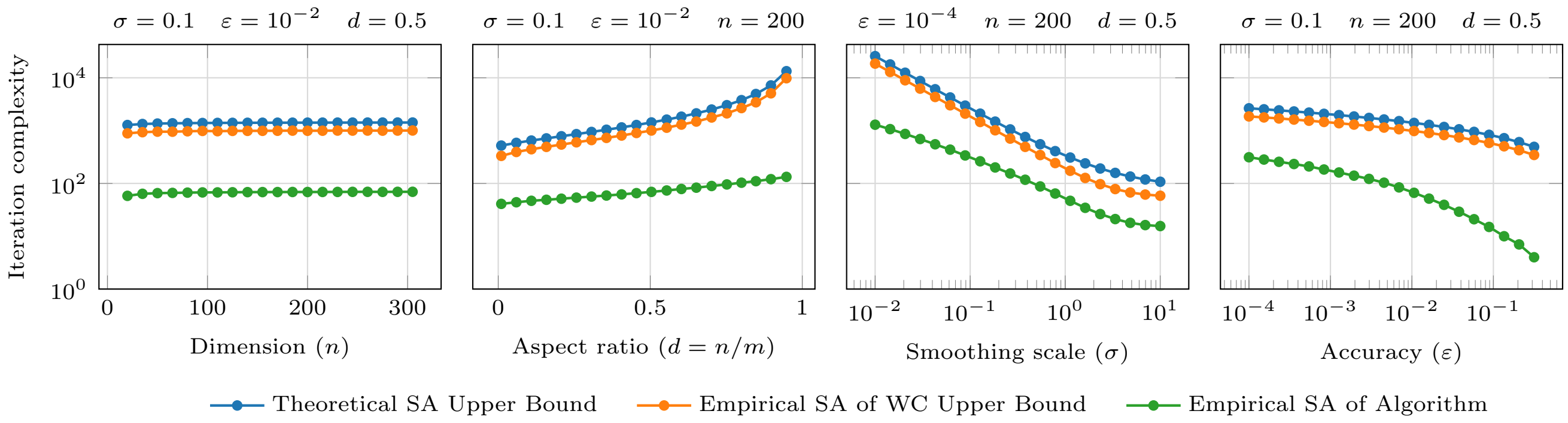}
    \caption{Smoothed complexity of Algorithm \eqref{alg:NAG} under Wishart perturbations.}
    \label{fig:NAG}
\end{figure}

\textbf{Numerical Illustration.} We use the same weighted Huber objective as for Gradient Descent (Section \ref{sec:GD1L}),
\[
    f(x)=\sum_{i=1}^n 2^{10(i-n/2)/n}\phi(x_i),
    \qquad
    \phi(t)=\begin{cases}
        \frac12t^2, & |t|\le 1,\\
        |t|-\frac12, & |t|>1.
    \end{cases}
\]
The results are shown in Figure \ref{fig:NAG}. We observe that the theoretical bound and the empirical expectation of the worst-case bound exhibit the same dependencies, up to constants. 

\begin{figure}[H]
    \centering
    \def\SAMethod{NAGQuadratic}
    \input{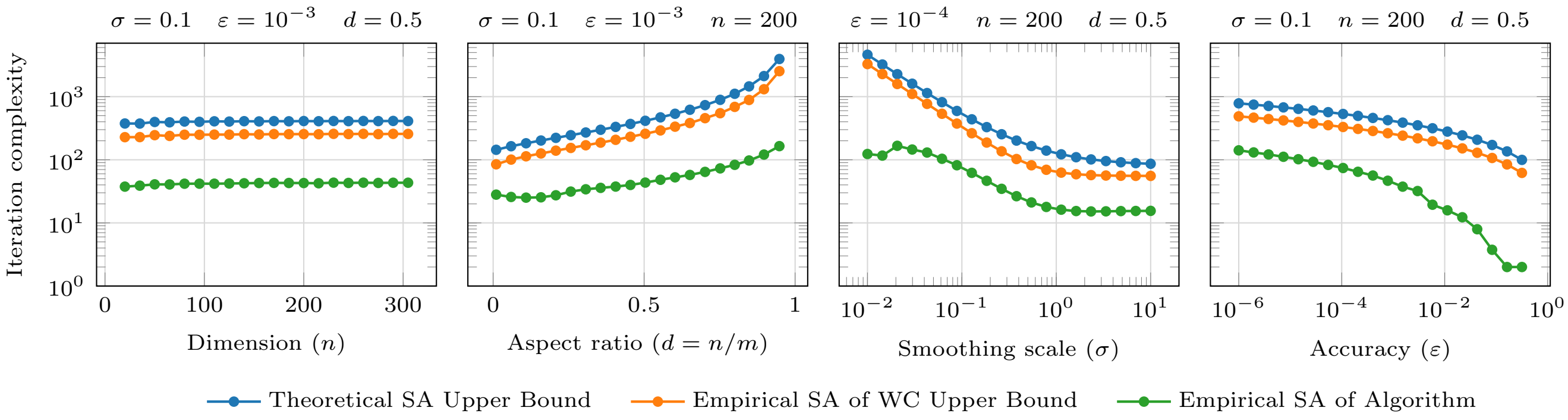}
    \caption{Smoothed complexity of Algorithm \eqref{alg:NAG} under Wishart perturbations.}
    \label{fig:NAGQuadratic}
\end{figure}

We additionally consider the degenerate ill-conditioned quadratic from Section \ref{sec:GD1L}, given by 
\[
    f(x)=\frac12\sum_{i=1}^{\lfloor n/2\rfloor}a_i x_i^2,
    \qquad a_i=2^{10(i-n/2)/n}.
\]
The results are shown in Figure \ref{fig:NAGQuadratic}. We note that our results correctly predict the dependencies in the shown variables, up to constant factors.

\subsection{Conjugate Gradient}\label{sec:CG}

We next turn to the \textit{Conjugate Gradient} method \citep{hestenes1952methods}, which exploits the quadratic structure of the objective. Let $\mathcal Q_\mu^L\subset \mathcal F_\mu^L$ denote the class of $\mu$-strongly convex quadratic functions $h(x)=\frac12x^TAx-b^Tx+c$ with $0\preceq \mu I\preceq A\preceq LI$. For $h\in \mathcal Q_\mu^L$, Conjugate Gradient generates, where $r_k=b-Ax_k$ and $p_0=r_0$, for $x_0\in \R^n$ and $k\ge 0$,
\begin{equation}\tag{CG}\label{alg:CG}
    \alpha_k=\frac{r_k^Tr_k}{p_k^TAp_k}, \quad x_{k+1}=x_k+\alpha_kp_k, \quad r_{k+1}=r_k-\alpha_kAp_k, \quad p_{k+1}=r_{k+1}+\frac{r_{k+1}^Tr_{k+1}}{r_k^Tr_k}p_k,
\end{equation}
where we terminate on $r_k=0$. 

As previously, we measure complexity in terms of relative function value accuracy, namely, for $\varepsilon\in (0, 1)$,
\[
    \mathbb K^\varepsilon_{\text{CG}}(h)=\min_{k\ge 0}\left\{k\colon \frac{h(x_k)-\inf h}{h(x_0)-\inf h}\le \varepsilon\right\},
\]
where $\mathbb K^\varepsilon_{\text{CG}}(h)=0$ if $h(x_0)=\inf h$.

We first recall a classic worst-case complexity result (see, for instance, \cite[Theorem 38.5]{trefethen_numerical_1997}). The rate corresponds to an accelerated rate, as presented in Section \ref{sec:NAG}, in the sense that it depends on $\sqrt{\kappa}$ rather than $\kappa$.

\begin{theorem}[Worst-Case Complexity of \eqref{alg:CG} - {\cite{trefethen_numerical_1997}}]
    Let $h\in \mathcal Q_\mu^L$ with $0<\mu\le L<+\infty$, and let $(x_k)$ be the iterates of \eqref{alg:CG}. Then, for $\kappa=\frac{L}{\mu}$,
    \[
        \frac{h(x_k)-\inf h}{h(x_0)-\inf h}=\frac{\|x_k-x_*\|^2_{\nabla^2 h}}{\|x_0-x_*\|^2_{\nabla^2 h}}\le 4\left(\frac{\sqrt{\kappa}-1}{\sqrt{\kappa}+1}\right)^{2k}.
    \]
    Therefore, the relative accuracy complexity is bounded by
    \[
        \sup_{h\in \mathcal Q_\mu^L}\K_{\text{CG}}^\varepsilon(h)\le \left\lceil\frac{\sqrt{\kappa}}{2}\log\frac4\varepsilon\right\rceil.
    \]
\end{theorem}

We now provide a smoothed complexity bound for this method.


\begin{theorem}[Smoothed Complexity of \eqref{alg:CG}]\label{thm:SA_CG}
    Let $0<\sigma, L<+\infty$ and $\varepsilon\in (0, 1)$. We denote by $f_W$ a perturbation of $f$, parametrized by $\sigma$, as defined in Equation \eqref{eq:to_be_computed}. Then, as $n\to \infty$, the smoothed complexity of Algorithm \eqref{alg:CG} over the function class $\mathcal Q_0^{L}$ is bounded by
    \[
        \sup_{f\in \mathcal Q_0^{L}}\ExpD{W}{\K^\varepsilon_{\text{CG}}(f_W)}\le 1+\mathcal O_d\left(\left(1+\sigma^{-1}\sqrt{L}\right)\log\frac4\varepsilon\right).
    \]
\end{theorem}
\begin{proof}
    Fix any $f\in \mathcal Q_0^L$. As for Gradient Descent, it holds that 
    \[
        \kappa_\sigma=\frac{Lm+\sigma^2s_{\max}(G)^2}{\sigma^2s_{\min}(G)^2}.
    \]
    Using the fact that $(a+b)^{1/2}\le a^{1/2}+b^{1/2}$, Cauchy-Schwarz and Theorem \ref{thm:moments}, we obtain
    \[
        \ExpD{W}{\K_{\text{CG}}^\varepsilon(f_W)}\le 1+\ExpD{W}{\left(\frac{Lm}{\sigma^2 s_{\min}(G)^2}+\frac{s_{\max}(G)^2}{s_{\min}(G)^2}\right)^{1/2}}\frac{\log(4/\varepsilon)}2\le 1+\mathcal O_d\left((1+\sigma^{-1}\sqrt L)\log\frac4\varepsilon\right). \qedhere
    \]
\end{proof}

\begin{remark}
    Smoothed analysis of Conjugate Gradient was previously studied in \cite{menon_smoothed_2016}, who study linear systems under random matrix perturbations, or equivalently, quadratic objectives under random quadratic perturbations. Their analysis heavily relies on spectral bounds, whereas ours is derived from moment bounds. This highlights a difference in the motivations: we aim to derive simple bounds directly from moment bounds, applicable across different optimization papers in a unified manner. Although the random matrix ensemble differs, we highlight that our complexity estimate is dimension-free, unlike theirs. 
\end{remark}

\textbf{Discussion.} We again observe a logarithmic dependency in $\varepsilon$, arising due to the perturbation, and a dependency in $\sigma^{-1}\sqrt{L}$, corresponding to the square root of the condition number.

\begin{figure}[H]
    \centering
    \def\SAMethod{CG}
    \input{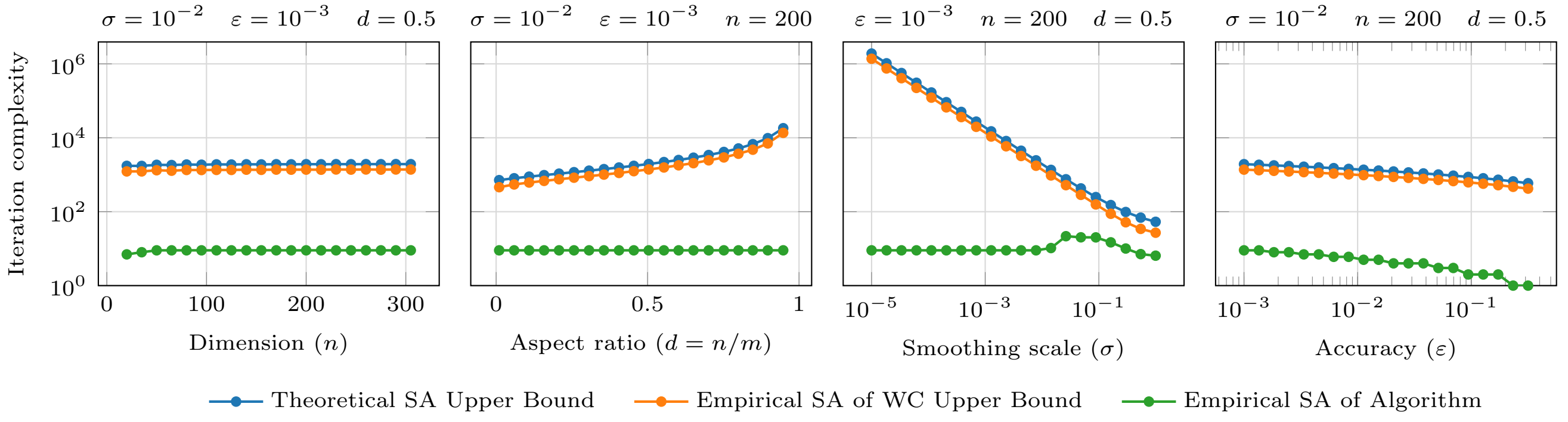}
    \caption{Smoothed complexity of Algorithm \eqref{alg:CG} under Wishart perturbations.}
    \label{fig:CG}
\end{figure}

\textbf{Numerical Illustration.} We illustrate Theorem \ref{thm:SA_CG} on the degenerate convex quadratic
\[
    f(x)=\frac12\sum_{i=1}^{\lfloor n/2\rfloor}a_i x_i^2,
    \qquad a_i=2^{10(i-n/2)/n}.
\]
The positive eigenvalues range from $2^{10(1-n/2)/n}$ to $L=2^{10(\lfloor n/2\rfloor-n/2)/n}$, while the remaining eigenvalues are zero.
The results are shown in Figure \ref{fig:CG}. The dotted line without markers shows CG applied directly to the unperturbed quadratic. This baseline is independent of $\sigma$ and $d$. Similarly to the previous section, we observe that the smoothed complexity bound is close to the empirical smoothed complexity of the worst-case bound, but very pessimistic compared to the empirical smoothed complexity of the run algorithm.


\subsection{Cubic Regularized Newton Method}\label{sec:CNM}

We next turn to the \textit{Cubic Regularized Newton Method} (CNM) \citep{nesterov_cubic_2006} on the class $\mathcal F_{\mu}^{L_H}(\Delta)$ of $\mu$-strongly convex functions with $L_H$-Lipschitz continuous Hessian and initialization bounded by $\Delta$, with the convention that $\mu=0$ corresponds to the merely convex bounded-below setting. For $h\in \mathcal F_{\mu}^{L_H}(\Delta)$, we consider, given $x_0\in \R^n$ and $k\ge 0$,
\begin{equation}\tag{CNM}\label{alg:CNM}
    s_k\in \argmin\left\{m_{x_k}(s)=h(x_k)+\nabla h(x_k)^Ts+\frac12s^T\nabla^2 h(x_k)s+\frac{L_H}{6}\|s\|^3\right\}, \quad x_{k+1}=x_k+s_k.
\end{equation}
We now measure complexity in terms of absolute function value accuracy, namely, for $\varepsilon>0$,
\[
    \mathbb K^\varepsilon_{\text{CNM}}(h)=\min_{k\ge 0}\left\{k\colon h(x_k)-\inf h\le \varepsilon\right\}.
\]
We first recall a worst-case complexity result, combining the slow phase from \cite[Theorem 5]{nesterov_cubic_2006} with the fast phase from \cite[Section 6]{nesterov_accelerating_2008}. This result is of the classical form $C+\log\log(1/\varepsilon)$.

\begin{theorem}[Worst-Case Complexity for \eqref{alg:CNM} - {\cite{nesterov_cubic_2006,nesterov_accelerating_2008}}]\label{thm:WC_CNM}
    Let $h\in \mathcal F_\mu^{L_H}(\Delta)$ for $0<\mu, L_H<+\infty$ and $0\le \Delta<+\infty$, and let $(x_k)$ be the iterates of \eqref{alg:CNM}. 
    Then, the absolute accuracy complexity is bounded by, for $\varepsilon>0$,
    \[
    \sup_{h\in \mathcal F_\mu^{L_H}(\Delta)}\K_{\text{CNM}}^\varepsilon(h)\le \K_{\slow}(\mu, \Delta)+\K_{\fast}^\varepsilon(\mu),
    \]
    where
    \[
        \K_{\slow}(\mu, \Delta)
        =\left\lceil 6.25\cdot2^{3/4}\cdot \left(\frac{\Delta L_H^2}{\mu^3}\right)^{1/4}\right\rceil,
        \quad \K_{\fast}^\varepsilon(\mu)=\left\lceil\log_{2,+}\log_{2,+}\left(\frac{\mu^3}{2L_H^2\varepsilon}\right)\right\rceil.
    \]
\end{theorem}

Stronger results, although with more complicated dependencies on the variables and more delicate assumptions, exist, and we purposefully do not make use of them here due to the additional complication. We refer the reader to \cite{cartis_adaptive_2011} for an initial overview, and to \cite{cartis_evaluation_2022} for a comprehensive survey.

As before, the worst-case bound is unbounded in the merely convex setting, whereas the random quadratic regularization makes the bound finite in expectation. This yields the following smoothed complexity guarantee.

\begin{theorem}[Smoothed Complexity for \eqref{alg:CNM}]\label{thm:SA_CNM}
    Let $0\le \Delta<+\infty$ and $0<\sigma,\varepsilon, L_H<+\infty$. We denote by $f_W$ a perturbation of $f$, parametrized by $\sigma$, as defined in Equation \eqref{eq:to_be_computed}. Then, as $n\to \infty$, the smoothed complexity of Algorithm \eqref{alg:CNM} over the function class $\mathcal F_{0}^{L_H}(\Delta)$ is bounded by
    \[
        \sup_{f\in \mathcal F_{0}^{L_H}(\Delta)}\ExpD{W}{\K^\varepsilon_{\text{CNM}}(f_W)}\le 2+\mathcal O_d\left(  L_H^{1/2}\left(\sigma^{-3/2}\Delta^{1/4}+D^{1/2}\sigma^{-1}\right)+\log_{2}\left(1+\log_2 \left(1+\frac{\sigma^6}{L_H^2\varepsilon}\right)\right)\right),
    \]
    where $\|x_0-\tilde x\|\le D$.
\end{theorem}
\begin{proof}
    Fix any $f\in \mathcal F_{0}^{L_H}(\Delta)$.
    We consider the fast and slow phases separately:
    \begin{itemize}
        \item \textbf{Slow Phase.} Using $(a+b)^{1/4}\le a^{1/4}+b^{1/4}$, Equations \eqref{eq:mu_L} and \eqref{eq:Del0}, Cauchy-Schwarz and Theorem \ref{coro:mixed_moments}, 
        {\small\begin{align*}
            \Exp{\K_{\slow}(\mu_\sigma, \Delta_\sigma)}
            &\le 1+\Exp{6.25\cdot2^{3/4}\cdot \left(\Delta+{\sigma^2D^2}s_{\max}(G)^2/(2m)\right)^{1/4}L_H^{1/2}\sigma^{-3/2}s_{\min}(G)^{-3/2}m^{3/4}} \\
            &\le 1+\mathcal O\left(\Delta^{1/4}L_H^{1/2}\sigma^{-3/2}\Exp{m^{3/4}s_{\min}(G)^{-3/2}}+\sigma^{-1}D^{1/2}L_H^{1/2}\Exp{\frac{s_{\max}(G)^{1/2}}{m^{1/4}}\frac{s_{\min}(G)^{-3/2}}{m^{-3/4}}}\right) \\
            &\le 1+\mathcal O_d(\Delta^{1/4}L_H^{1/2}\sigma^{-3/2}+\sigma^{-1}D^{1/2}L_H^{1/2}).
        \end{align*}}

        \item \textbf{Fast Phase.} By Jensen's inequality, Equation \eqref{eq:mu_L} and the fact that $s_{\min}(G)\le s_{\max}(G)$,
        \[
            \Exp{\K_\fast^\varepsilon(\mu_\sigma)}\le 1+\Exp{\log_2\left(1+\log_2\left(1+\frac{\sigma^6}{2L_H^2\varepsilon}\frac{s_{\min}(G)^6}{m^3}\right)\right)}\le 1+\log_2\left(1+\log_2\left(1+\frac{\sigma^6}{2L_H^2\varepsilon}\Exp{\frac{s_{\max}(G)^6}{m^3}}\right)\right).
        \]
        We conclude by Theorem \ref{coro:mixed_moments}.
        \qedhere
    \end{itemize}
\end{proof}

\textbf{Discussion.} We observe that the smoothed complexity bound takes the form $C_\sigma+\log\log(1/\varepsilon)$, where the constant $C_\sigma$ absorbs $L_H, \Delta, D, \sigma$, with $C_\sigma$ acting like $\sigma^{-3/2}$. The quadratic phase is only valid for $\varepsilon$ sufficiently small.

\begin{figure}[H]
    \centering
    \def\SAMethod{CNMLogistic}
    \input{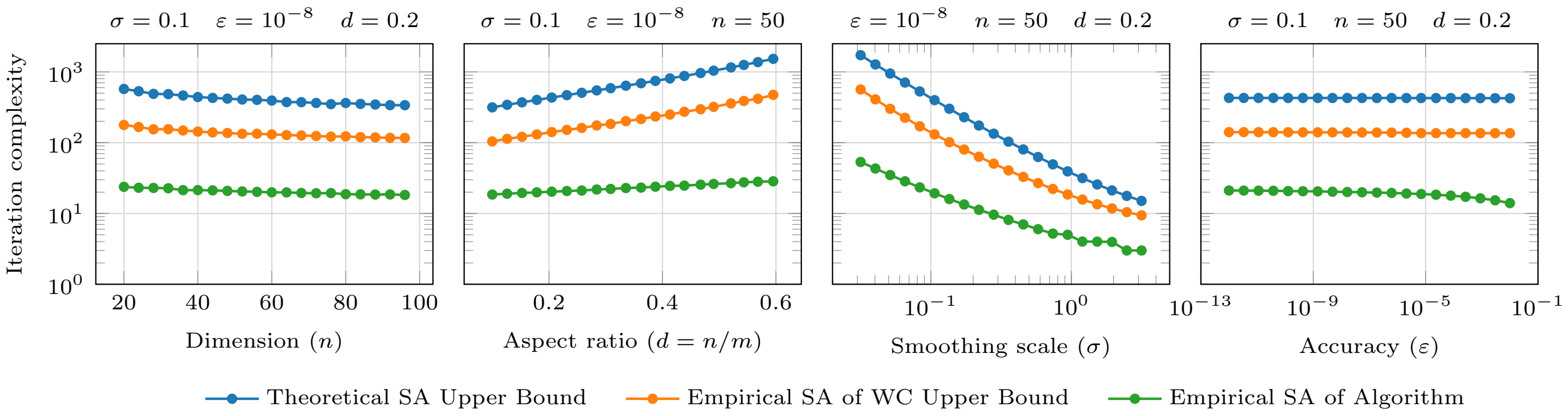}
    \caption{Smoothed complexity of Algorithm \eqref{alg:CNM} under Wishart perturbations.}
    \label{fig:CNM}
\end{figure}

\textbf{Numerical Illustration.} We illustrate Theorem \ref{thm:SA_CNM} on the logistic regression objective, for $r\ge 2$,
\[
    f(x)=\frac{1}{r}\sum_{i=1}^r\left(\log(1+e^{b_i^Tx})-\frac34b_i^Tx\right).
\]
This is the mean binary logistic loss for $4r$ observations, where each feature vector $b_i$ occurs four times, with three labels equal to one and one equal to zero. We define $B=(b_1, \ldots, b_r)^T=\text{diag}(s_i)Q^T$, where $Q\in \R^{n\times r}$ has orthonormal columns obtained through a QR factorization of a standard Gaussian matrix, and $s_i=2^{\frac{3(i-1)}{r-1}-2}$. A minimizer and minimal value is given by 
\[
    x_*=Q\text{diag}(s_i^{-1})\log 3\cdot \mathbf 1\quad \text{and}\quad \min f=-\frac34\log \frac34-\frac14\log \frac14.
\]
We easily see $\nabla^2 f$ is positive semidefinite and rank-deficient, making $f$ convex but not strongly convex. Moreover, we note that $f$ is $L_G$-smooth and has $L_H$-Lipschitz Hessian, for 
\[
    L_G=\frac{1}{r}\quad \text{and}\quad L_H=\frac{4}{3\sqrt{3}r}.
\]
The results are shown in Figure \ref{fig:CNM}. We note that, unlike in previous experiments, the Lipschitz constants $L_G$ and $L_H$ depend on the dimension of the problem $n=4r$, and therefore the curves are not constant in $n$, although the upper bounds are dimension-independent. The dependency in $d$ is pessimistic but with the correct trends, and the dependency in $\sigma$ matches up to constants.

\subsection{Cubic Regularized Newton Method under Radial Monotonicity}\label{sec:CNM_PL}

We now consider Algorithm \eqref{alg:CNM} beyond convexity, as in Section \ref{sec:GD1L_PL}. We denote by $\mathcal F_{\text{PL}, \eta}^{L_H}(\Delta)$ the class of twice continuously differentiable functions with $L_H$-Lipschitz continuous Hessian, bounded below, satisfying gradient domination with parameter $\eta$, and with initialization bounded by $h(x_0)-\min h\le \Delta$. We use the absolute function value accuracy complexity $\mathbb K_{\text{CNM}}^\varepsilon(h)$, as in Section \ref{sec:CNM}. As for Gradient Descent, gradient domination is sufficient to retain linear convergence beyond convexity. This is formalized through the result of \cite[Theorem 7]{nesterov_cubic_2006}.

\begin{theorem}[{Worst-Case Complexity of \eqref{alg:CNM} under Gradient Domination - \cite{nesterov_cubic_2006}}]\label{thm:WC_CNM_PL}
    Let $h\in \mathcal F_{\text{PL}, \eta}^{L_H}(\Delta)$ with $0<\eta, L_H<+\infty$ and $0\le \Delta<+\infty$, and let $(x_k)$ be the iterates of \eqref{alg:CNM}. Then, the absolute function value accuracy complexity, for every $\varepsilon>0$, is bounded by 
    \[
        \sup_{h\in \mathcal F_{\text{PL}, \eta}^{L_H}(\Delta)}\mathbb K^\varepsilon_{\text{CNM}}(h)\le \mathbb K_{\text{slow}}(\eta, \Delta)+\mathbb K^\varepsilon_{\text{fast}}(\eta),
    \]
    where
    \[
        \mathbb K_{\text{slow}}(\eta, \Delta)= \left\lceil8\left(\frac{\Delta}{\omega}\right)^{1/4}\right\rceil, \quad \mathbb K^\varepsilon_{\text{fast}}(\eta)=\left\lceil\log_{4/3,+}\log_{2,+}\left(\frac{\omega}{\varepsilon}\right)\right\rceil\quad\text{and}\quad \omega=\frac{\eta^3}{2592L_H^2}.
    \]
\end{theorem}

We denote by $\mathcal F_{\text{rad}}^{L_G, L_H}(\Delta, D; \tilde x)$ the subclass of $\mathcal F_{\text{rad}}^{L_G}(\tilde x)$ from Section \ref{sec:GD1L_PL} that is also twice continuously differentiable with $L_H$-Lipschitz Hessian and bounded initialization $f(x_0)-\min f\le \Delta$ and $\|x_0-\tilde x\|\le D$. The quadratic perturbation preserves the Lischitz constant of the Hessian $L_H$, and the arguments in Section \ref{sec:GD1L_PL} yield
\[
    \mu_\sigma=\sigma^2\lambda_{\min}(W), \quad L_{G, \sigma}=L_G+\sigma^2\lambda_{\max}(W), \quad \eta_\sigma=\frac{\mu_\sigma^2}{L_{G, \sigma}}.
\]
In particular, $f_W\in \mathcal F_{\text{PL}, \eta_\sigma}^{L_{H}}(\Delta_\sigma)$. This allows us to establish the following smoothed complexity guarantee. 

\begin{theorem}[Smoothed Complexity of \eqref{alg:CNM} under Radial Monotonicity]\label{thm:SA_CNM_PL}
    Let $0\le \Delta, D<+\infty$, $0<\sigma, L_G, L_H, \varepsilon<+\infty$ and $\tilde x\in \R^n$. We denote by $f_W$ a perturbation of $f$, parametrized by $\sigma$, as defined in Equation \eqref{eq:to_be_computed}. Then, as $n\to \infty$, the smoothed complexity of Algorithm \eqref{alg:CNM} over the function class $\mathcal F_{\text{rad}}^{L_G, L_H}(\Delta, D; \tilde x)$ is bounded by 
    \begin{align*}
        \sup_{f\in \mathcal F_{\text{rad}}^{L_G, L_H}(\Delta, D; \tilde x)}\ExpD{W}{\mathbb K^\varepsilon_{\text{CNM}}(f_W)}\le 2+\mathcal O_d& \left(L_H^{1/2}(1+\sigma^{-2}L_G)^{3/4}(\sigma^{-3/2}\Delta^{1/4}+\sigma^{-1}D^{1/2})\right. \\
        &\left.\quad +\log_{4/3}\left(1+\log_2\left(1+\frac{\sigma^6}{L_H^2\varepsilon}\right)\right)\right).
    \end{align*}
\end{theorem}
\begin{proof}
    Fix any $f\in \mathcal F_{\text{rad}}^{L_G, L_H}(\Delta, D; \tilde x)$. We consider the slow and the fast phase separately.
    \begin{itemize}
        \item \textbf{Slow Phase.} By Theorem \ref{thm:WC_CNM_PL}, it holds that
        \[
            \Exp{\mathbb K_{\text{slow}}(\eta_\sigma, \Delta_\sigma)}\le 1+\ExpD{W}{8\left(\frac{\Delta_\sigma}{\omega_\sigma}\right)^{1/4}}=1+\mathcal O_d\left(\ExpD{W}{\Delta_\sigma^{1/4}\eta_\sigma^{-3/4}L_H^{1/2}}\right)=1+\mathcal O_d\left(\ExpD{W}{\Delta_\sigma^{1/4}\mu_\sigma^{-3/2}L_{G,\sigma}^{3/4}L_H^{1/2}}\right).
        \]
        Using subadditivity of powers $<1$, we get
        \begin{align*}
            \Exp{\mathbb K_{\text{slow}}(\eta_\sigma, \Delta_\sigma)}
            &\le1+ \mathcal O_d\left(\ExpD{W}{(\Delta+\sigma^2D^2\lambda_{\max}(W))^{1/4}\sigma^{-3}\lambda_{\min}(W)^{-3/2}(L_G+\sigma^2\lambda_{\max}(W))^{3/4}L_H^{1/2}}\right) \\
            &\le1+ \mathcal O_d\left(\ExpD{W}{\frac{\sigma^{-3/2}\Delta^{1/4}+\sigma^{-1}D^{1/2}\lambda_{\max}(W)^{1/4}}{\lambda_{\min}(W)^{3/4}}\frac{(\sigma^{-2}L_G+\lambda_{\max}(W))^{3/4}}{\lambda_{\min}(W)^{3/4}}L_H^{1/2}}\right).
        \end{align*}
        We conclude by Cauchy-Schwarz and Theorem \ref{thm:moments}.
        \item \textbf{Fast Phase.} By Jensen's inequality, we get 
        \[
            \Exp{\mathbb K^\varepsilon_{\text{fast}}}(\eta_\sigma)\le 1+\Exp{\log_{4/3}\left(1+\log_2\left(1+\frac{\omega_\sigma}{\varepsilon}\right)\right)}\le 1+\log_{4/3}\left(1+\log_2\left(1+\frac{\Exp{\omega_\sigma}}{\varepsilon}\right)\right).
        \]
        Since
        \[
            \Exp{\omega_\sigma}=\mathcal O_d\left(\Exp{\frac{\eta_\sigma^3}{L_H^2}}\right)=\mathcal O_d\left(L_H^{-2}\Exp{\mu_\sigma^6L_{G, \sigma}^{-3}}\right)=\mathcal O_d\left(L_H^{-2}\Exp{\frac{\sigma^{12}\lambda_{\min}(W)^6}{(L_G+\sigma^2\lambda_{\max}(W))^3}}\right),
        \]
        which implies
        \[
            \Exp{\omega_\sigma}\le \mathcal O_d\left(L_H^{-2}\sigma^{6}\Exp{\frac{\lambda_{\min}(W)^6}{\lambda_{\max}(W)^3}}\right)\le \mathcal O_d\left(L_H^{-2}\sigma^{6}\Exp{\lambda_{\max}(W)^3}\right)=\mathcal O_d\left(\frac{\sigma^6}{L_H^2}\right),
        \]
        we get the wanted bound. \qedhere
    \end{itemize}
\end{proof}

\textbf{Discussion.} Theorem \ref{thm:SA_CNM_PL} is very closely related to Theorem \ref{thm:SA_CNM}, whose upper bound only differs in a factor $(1+\sigma^{-2}L_G)^{3/4}$. This is comparable to the PL extension of Theorem \ref{thm:SA_GD1L}, provided in Theorem \ref{thm:SA_GD1L_PL}, which also only differed in one additional factor. 

\begin{figure}[H]
    \centering
    \def\SAMethod{CNMRadial}
    \input{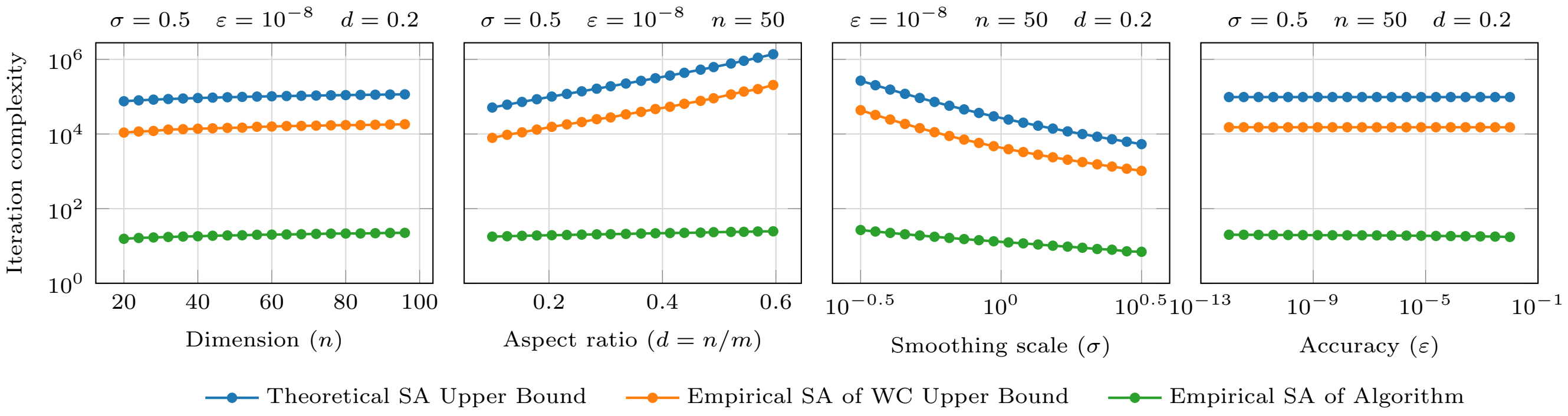}
    \caption{Smoothed complexity of Algorithm \eqref{alg:CNM} under Wishart perturbations.}
    \label{fig:CNM_PL}
\end{figure}

\textbf{Numerical Illustration.} We consider the same radially monotone function as in Section \ref{sec:GD1L_PL}, namely 
\[
    f(x)=\sum_{i=1}^n\left(1-e^{-x_i^2/2}\right), \quad \tilde x=0.
\]
The results are shown in Figure \ref{fig:CNM_PL}. As in Section \ref{sec:GD1L_PL}, we observe that the PL extension of our result is less tight in the dependencies than the convex counterpart. However, we still observe similar dependencies in all plotted variables. We note that the theoretical upper bounds are worse than the ones presented in Section \ref{sec:GD1L_PL} by a constant factor.

\subsection{Hybrid Gradient--Cubic Newton Method}\label{sec:Hyb}

We combine Gradient Descent and the Cubic Regularized Newton Method on the class $\mathcal F_{\mu}^{L_G,L_H}(\Delta)$ of bounded-below $\mu$-strongly convex, $L_G$-smooth functions with $L_H$-Lipschitz continuous Hessian and initializations bounded by $h(x_0)-\inf h\le \Delta$. For $h\in \mathcal F_{\mu}^{L_G,L_H}(\Delta)$, we run Algorithm \eqref{alg:GD1L} with step-size $1/L_G$ until
\[
    \|\nabla h(x_k)\|^2\le \frac{\mu^4}{2L_H^2},
\]
after which we run Algorithm \eqref{alg:CNM}. By strong convexity, the switching criterion implies $h(x_k)-\min h\le \frac{\mu^3}{4L_H^2}$, so that the method enters the fast phase of Theorem \ref{thm:WC_CNM}. We still measure complexity in terms of absolute function value accuracy, namely
\[
    \mathbb K^\varepsilon_{\text{Hyb}}(h)=\min_{k\ge 0}\left\{k\colon h(x_k)-\inf h\le \varepsilon\right\}.
\]
This improves the rate of the previous section to a rate of the form $\log(C)+\log\log(1/\varepsilon)$.

\begin{theorem}[Worst-Case Complexity of Hybrid Method]\label{thm:WC_Hyb}
    Let $h\in \mathcal F_{\mu}^{L_G,L_H}(\Delta)$ with $0<\mu\le L_G<+\infty$, $0<L_H<+\infty$, $0\le \Delta<+\infty$. 
    Then, the absolute accuracy complexity is bounded by, for $\varepsilon>0$,
    \[
    \sup_{h\in \mathcal F_\mu^{L_G,L_H}(\Delta)}\K_{\text{Hyb}}^\varepsilon(h)\le \K_{\slow}(\mu,L_G,\Delta)+\K_{\fast}^\varepsilon(\mu),
    \]
    where
    \[
        \K_{\slow}(\mu,L_G,\Delta)
        =\left\lceil\frac{L_G}{\mu}\log_+ \frac{4L_G L_H^2\Delta}{\mu^4}\right\rceil,
        \quad \K_{\fast}^\varepsilon(\mu)=\left\lceil\log_{2,+}\log_{2,+}\frac{\mu^3}{2L_H^2\varepsilon}\right\rceil.
    \]
\end{theorem}
\begin{proof}
    Let $\kappa=L_G/\mu$. By $L_G$-smoothness, $\|\nabla h(x_k)\|^2\le 2L_G(h(x_k)-\min h)$, so the switching criterion is guaranteed once
    \[
        \frac{h(x_k)-\inf h}{h(x_0)-\inf h}\le \frac{\mu^3}{4\kappa L_H^2\Delta}.
    \]
    Theorem \ref{thm:WC_GD1L} shows that this takes at most
    \[
        \K_{\slow}(\mu,L_G,\Delta)=\left\lceil\kappa\log_+\frac{4\kappa L_H^2\Delta}{\mu^3}\right\rceil=\left\lceil\frac{L_G}{\mu}\log_+\frac{4L_GL_H^2\Delta}{\mu^4}\right\rceil
    \]
    iterations. Once the switching criterion holds, strong convexity gives $h(x_k)-\min h\le \frac{\mu^3}{4L_H^2}$. The fast phase therefore follows from Theorem \ref{thm:WC_CNM}. 
\end{proof}

\begin{theorem}[Smoothed Complexity of Hybrid Method]\label{thm:SA_Hyb}
    Let $0\le \Delta<+\infty$ and $0<\sigma,\varepsilon,L_G,L_H<+\infty$. We denote by $f_W$ a perturbation of $f$, parametrized by $\sigma$, as defined in Equation \eqref{eq:to_be_computed}. Then, as $n\to \infty$, the smoothed complexity of the hybrid algorithm over the function class $\mathcal F_{0}^{L_G,L_H}(\Delta, D)$ is bounded by
    \begin{align*}
        \sup_{f\in \mathcal F_{0}^{L_G,L_H}(\Delta)}\ExpD{W}{\K_{\text{Hyb}}^\varepsilon(f_W)}
        \le 2+&\mathcal O_d\Bigl((1+\sigma^{-2}L_G)\log(2+L_H^2(1+\sigma^{-2}L_G)(\Delta\sigma^{-6}+D^2\sigma^{-4})) \\ &\qquad +\left. \log_{2}\left(1+\log_2\left(1+ \frac{\sigma^6}{L_H^2\varepsilon}\right)\right)\right),
    \end{align*}
    where $\|x_0-\tilde x\|\le D$.
\end{theorem}
\begin{proof}
    Fix any $f\in \mathcal F_{0}^{L_G,L_H}(\Delta,D)$.
    We consider the slow and fast phases separately:
    \begin{itemize}
        \item \textbf{Slow Phase.} By Theorem \ref{thm:WC_Hyb} and Equations \eqref{eq:mu_L} and \eqref{eq:Del0},
        \begin{align*}
            \Exp{\K_{\slow,\sigma}}
            &\le1+ \Exp{\kappa_\sigma\log\left(1+\frac{2\kappa_\sigma L_H^2m^2(2m\Delta+\sigma^2D^2s_{\max}(G)^2)}{\sigma^{6}s_{\min}(G)^{6}}\right)}
            =1+\Exp{\kappa_\sigma\log(1+2A_\sigma)},
        \end{align*}
        where
        \[
            A_\sigma=\frac{\kappa_\sigma L_H^2m^2(2m\Delta+\sigma^2D^2s_{\max}(G)^2)}{\sigma^{6}s_{\min}(G)^{6}}.
        \]
        Here we used $\log_+(x)\le \log(1+x)$. Since $\log(1+2x)\le \log(1+\alpha)+\frac{x}{\alpha}$ for any $\alpha>0$, setting $\alpha=\sqrt{\Exp{A_\sigma^2}}$ yields
        \[
            \Exp{\K_{\slow,\sigma}}\le 1+\Exp{\kappa_\sigma}\log(1+\sqrt{\Exp{A_\sigma^2}})+\frac{\Exp{\kappa_\sigma A_\sigma}}{\sqrt{\Exp{A_\sigma^2}}}.
        \]
        By Cauchy-Schwarz,
        \[
            \Exp{\K_{\slow,\sigma}}\le 1+\sqrt{\Exp{\kappa_\sigma^2}}\left(\log(1+\sqrt{\Exp{A_\sigma^2}})+1\right).
        \]
        It holds that 
        \[
            \ExpD{W}{\kappa_\sigma^2}
        \le 2L_G^2\sigma^{-4}\ExpD{W}{\frac{m^2}{s_{\min}(G)^4}}
        +2\ExpD{W}{\frac{s_{\max}(G)^4}{s_{\min}(G)^4}}
        \le\mathcal O_d(1+\sigma^{-4}L_G^2),
        \]
        and, since
        \begin{align*}
            A_\sigma
            &=L_H^2\frac{m^2(L_Gm+\sigma^2s_{\max}(G)^2)(2m\Delta+\sigma^2D^2s_{\max}(G)^2)}{\sigma^8s_{\min}(G)^8} \\
            &=L_H^2\left(
            \frac{2m^4L_G\Delta}{\sigma^8s_{\min}(G)^8}
            +\frac{2m^3\Delta s_{\max}(G)^2}{\sigma^6s_{\min}(G)^8}
            +\frac{L_Gm^3D^2s_{\max}(G)^2}{\sigma^6s_{\min}(G)^8}
            +\frac{m^2s_{\max}(G)^4D^2}{\sigma^4s_{\min}(G)^8}
            \right),
        \end{align*}
        it holds, by Cauchy-Schwarz and Theorem \ref{thm:moments},
        \begin{align*}
            \Exp{A_\sigma^2}
            &\le \mathcal O_d\left(L_H^4\left(L_G^2\Delta^2 \sigma^{-16}+\Delta^2 \sigma^{-12}+L_G^2D^4\sigma^{-12}+D^4\sigma^{-8}\right)\right) \\
            &= \mathcal O_d\left(L_H^4(1+\sigma^{-2}L_G)^2(\Delta^2\sigma^{-12}+D^4\sigma^{-8})\right).
        \end{align*}
        Therefore, 
        \[
            \Exp{\K_{\slow,\sigma}}\le 1+\mathcal O_d\left((1+\sigma^{-2}L_G)\left(\log(1+L_H^2(1+\sigma^{-2}L_G)(\Delta\sigma^{-6}+D^2\sigma^{-4}))+1\right)\right).
        \]
        \item \textbf{Fast Phase.} This follows exactly as in the proof of Theorem \ref{thm:SA_CNM}. \qedhere
    \end{itemize}
\end{proof}

\textbf{Discussion.} As expected, most of the terms that were part of the constant in Section \ref{sec:CNM} are now under the logarithm. Specifically, the only dependency that is not under a logarithm is $\sigma^{-2}L_G$, which acts as a condition number, and was also outside of the logarithm in Theorem \ref{thm:WC_CNM}. The quadratic convergence part remains the same, which is to be expected, as the algorithm is identical once the quadratic convergence region is attained.

\begin{figure}[H]
    \centering
    \def\SAMethod{Hybrid}
    \input{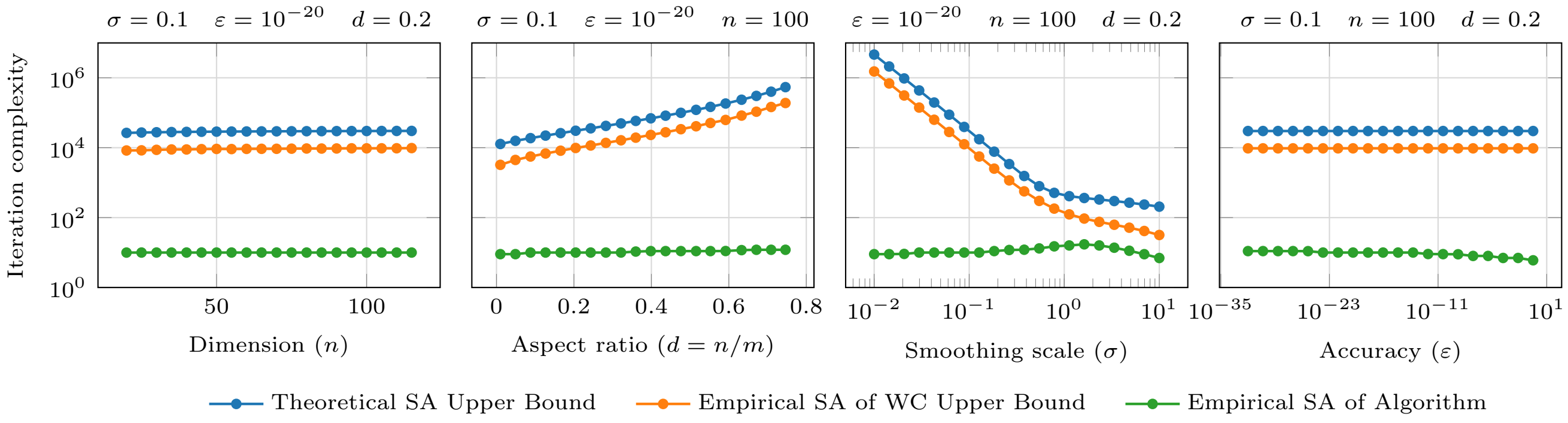}
    \caption{Smoothed complexity of the Hybrid Gradient--Cubic Newton method under Wishart perturbations.}
    \label{fig:Hybrid}
\end{figure}

\textbf{Numerical Illustration.} We illustrate Theorem \ref{thm:SA_Hyb} on the convex function
\[
    f(x)=\sum_{i=1}^n\left(\sqrt{1+x_i^2}-1\right),
\]
whose gradient is $1$-Lipschitz and whose Hessian is $1$-Lipschitz. The results are shown in Figure \ref{fig:Hybrid}. We observe that our bounds do not capture the tendency in the empirical iteration counters with $\sigma$ or $d$, but that they agree on the dimension-independence.


\subsection{Gauss-Newton Method - Convex Case}\label{sec:GNConvex}

We next consider \textit{Nonlinear Least Squares} problems $h(x)=\frac12\|r(x)\|^2$, and write $J(x)=\nabla r(x)$. Let $\mathcal F_{\mu,\underline\gamma}^{L,\bar\gamma}$ denote the class for which $f$ is $\mu$-strongly convex and $L$-smooth, while $\underline \gamma^2 I\preceq J(x)^TJ(x)\preceq \bar \gamma^2 I$ for every $x$. As before, $\mu=\underline\gamma=0$ denotes the merely convex, possibly rank-deficient setting. For $f\in\mathcal F_{\mu,\underline\gamma}^{L,\bar\gamma}$, we consider, for $x_0\in \R^n$ and $k\ge 0$,
\begin{equation}\tag{GN}\label{alg:GN}
    J_k=J(x_k),\quad g_k=\nabla f(x_k),\quad (J_k^TJ_k)s_k=-g_k,\quad x_{k+1}=x_k+\frac{\underline\gamma^2}{L}s_k.
\end{equation}
We measure complexity in terms of relative function value accuracy, namely, for $\varepsilon\in (0, 1)$,
\[
    \mathbb K^\varepsilon_{\text{GNConvex}}(h)=\min_{k\ge 0}\left\{k\colon \frac{h(x_k)-\inf h}{h(x_0)-\inf h}\le \varepsilon\right\},
\]
where $\mathbb K^\varepsilon_{\text{GNConvex}}(h)=0$ if $h(x_0)=\inf h$. Though we do not write it, for least squares problems with let the residual $r$ be part of the instance. 

We first recall the corresponding worst-case complexity. The proof, which follows the standard descent argument underlying Gauss--Newton analyses \citep{nocedal_numerical_2006,cartis_evaluation_2013,nesterov_modified_2003}, is included in Appendix \ref{app:proof_GN_Convex}.

\begin{theorem}[Worst-Case Complexity of \eqref{alg:GN}]\label{thm:GNConvex}
    Let $h\in \mathcal F_{\mu,\underline\gamma}^{L,\bar\gamma}$ for $0<\mu\le L<+\infty$, $0<\underline\gamma\le\bar\gamma<+\infty$, and let $(x_k)$ be the iterates of \eqref{alg:GN}. Then,
    \[
        \frac{h(x_k)-\inf h}{h(x_0)-\inf h}
        \le \left(1-\frac{\mu\underline\gamma^2}{L\bar\gamma^2}\right)^k.
    \]
    In particular, for $M=\max\{L,\bar\gamma^2\}$, $\ell=\min\{\mu,\underline\gamma^2\}$ and $R=M/\ell$, the relative accuracy complexity is bounded by, for $\varepsilon\in (0, 1)$,
    \[
        \sup_{h\in \mathcal F_{\mu,\underline\gamma}^{L,\bar\gamma}}\K_{\text{GNConvex}}^\varepsilon(h)
        \le \left\lceil\frac{L\bar\gamma^2}{\mu\underline\gamma^2}\log(1/\varepsilon)\right\rceil
        \le \left\lceil R^2\log(1/\varepsilon)\right\rceil.
    \]
\end{theorem}

In the merely convex setting, both the curvature of $f$ and the smallest singular value of $J$ may vanish. For the perturbed problem, the quadratic regularization can be written as an augmentation of the residual,
\begin{equation*}
    r_\sigma(x)=\begin{pmatrix}r(x)\\ \sigma W^{1/2}(x-\tilde x)\end{pmatrix},
    \qquad
    J_\sigma(x)=\begin{pmatrix}J(x)\\ \sigma W^{1/2}\end{pmatrix}.
\end{equation*}
Consequently,
\begin{equation}\label{eq:R_sigma}
    \frac{L_\sigma\bar\gamma_\sigma^2}{\mu_\sigma\underline\gamma_\sigma^2}
    \le \frac{(L+\sigma^2\lambda_{\max}(W))(\bar\gamma^2+\sigma^2\lambda_{\max}(W))}
    {\sigma^4\lambda_{\min}(W)^2}.
\end{equation}
This yields the following smoothed complexity guarantee.

\begin{theorem}[Smoothed Complexity of \eqref{alg:GN}]\label{thm:GNConvex_SA}
    Let $0<\sigma, L,\bar\gamma<+\infty$ and $\varepsilon\in(0,1)$. We denote by $f_W$ a perturbation of $f$, parametrized by $\sigma$, as defined in Equation \eqref{eq:to_be_computed}. Then, as $n\to\infty$, the smoothed complexity of the Algorithm \eqref{alg:GN} over the function class $\mathcal F_{0,0}^{L,\bar\gamma}$ is bounded by
    \[
        \sup_{f\in \mathcal F_{0,0}^{L,\bar\gamma}}\ExpD{W}{\K^\varepsilon_{\text{GNConvex}}(f_W)}
        \le1+ \mathcal O_d\left((1+\sigma^{-2}L)(1+\sigma^{-2}\bar\gamma^2)\log\frac1\varepsilon)\right).
    \]
\end{theorem}
\begin{proof}
    Fix any $f\in \mathcal F_{0,0}^{L,\bar\gamma}$.
    Theorem \ref{thm:GNConvex} and the perturbed parameters in Equation \eqref{eq:R_sigma} give
    \[
        \K_{\text{GNConvex}}^\varepsilon(f_W)
        \le 1+\frac{(L+\sigma^2\lambda_{\max}(W))(\bar\gamma^2+\sigma^2\lambda_{\max}(W))}
        {\sigma^4\lambda_{\min}(W)^2}\log\frac1\varepsilon.
    \]
    Expanding the numerator, applying Cauchy--Schwarz to the mixed moments and using Theorem \ref{thm:moments} gives the result. \qedhere
\end{proof}


\textbf{Discussion.} The quantities $\sigma^{-2}L$ and $\sigma^{-2}\bar\gamma^2$ both describe the conditioning scales. We therefore observe a dependency in the square condition number in the smoothed bound. The perturbation allows to recover linear convergence, meaning the upper bound scales logarithmically with $\varepsilon$, even though the degenerate setting only enjoys a sublinear convergence rate. 

\begin{figure}[H]
    \centering
    \def\SAMethod{GNConvex}
    \input{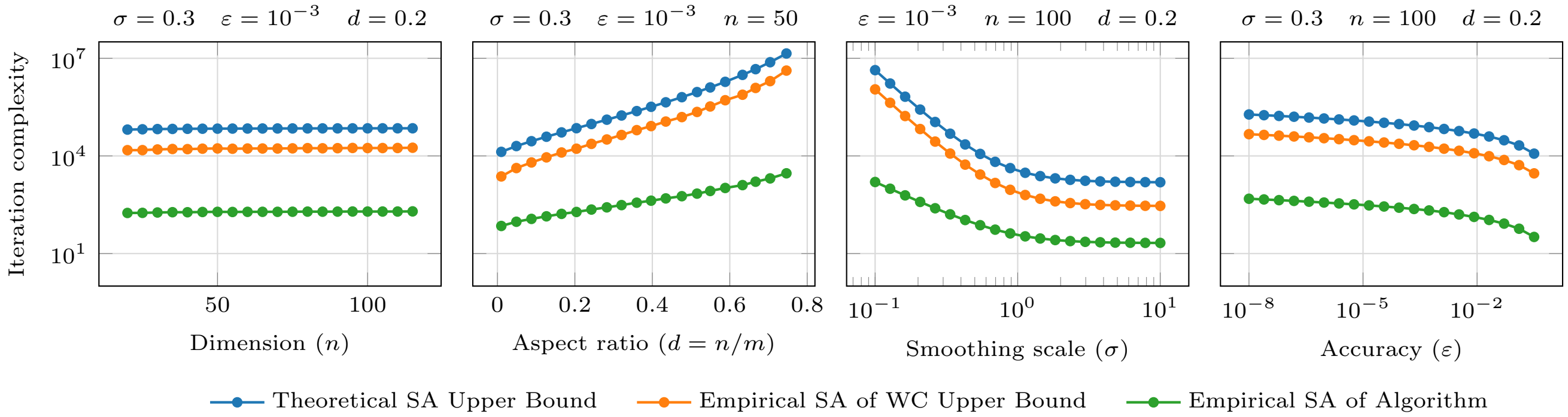}
    \caption{Smoothed complexity of Algorithm \eqref{alg:GN} under Wishart perturbations.}
    \label{fig:GNConvex}
\end{figure}

\textbf{Numerical Illustration.} We illustrate Theorem \ref{thm:GNConvex_SA} on the rank-one nonlinear least-squares objective
\[
    f(x)=\frac12\left(t(x)^2+\frac14\sin^2(t(x))\right),
    \qquad t(x)=\frac1{\sqrt n}\sum_{i=1}^n x_i.
\]
Indeed, $f$ is convex since its second derivative along $t$ is $1+\frac14\cos(2t)\ge\frac34$, while the residual Jacobian has rank one. 
The results are shown in Figure \ref{fig:GNConvex}. We observe that the bounds and the empirical iteration count increase with $d$ and exhibit a similar dependency in $\varepsilon$, up to constants. This again shows an agreement between the empirical and theoretical observations about the highly overparametrized perturbation model being favourable. However, the empirical iteration count decreases more slowly with $\sigma$ than the bounds for small perturbations, so that the theoretical dependency in $\sigma$ is pessimistic on this example. 

\subsection{Gauss-Newton Method - Convex Case with Armijo Line Search}\label{sec:GNConvexArmijo}

We consider the same class $\mathcal F_{\mu,\underline\gamma}^{L,\bar\gamma}$ as in Section \ref{sec:GNConvex}, replacing the fixed step-size in Algorithm \eqref{alg:GN} by an Armijo line search. Specifically, given the Armijo parameters $\beta,\tau\in(0,1)$, for $h\in\mathcal F_{\mu,\underline\gamma}^{L,\bar\gamma}$, consider, for $x_0\in\R^n$ and $k\ge 0$,
\begin{equation}\tag{GNArmijo}\label{alg:GNConvexArmijo}
    J_k=J(x_k),\quad g_k=\nabla h(x_k), \quad (J_k^TJ_k)s_k=-g_k,\qquad x_{k+1}=x_k+\alpha_ks_k,
\end{equation}
where $\alpha_k$ is the largest value in $1,\tau,\tau^2,\ldots$ satisfying
\begin{equation*}
    h(x_k+\alpha_ks_k)\le h(x_k)+\beta\alpha_kg_k^Ts_k.
\end{equation*}
We again measure complexity in terms of relative function value accuracy, namely, for $\varepsilon\in(0,1)$,
\[
    \mathbb K^\varepsilon_{\text{GNConvexArmijo}}(h)=\min_{k\ge0}\left\{k\colon\frac{h(x_k)-\inf h}{h(x_0)-\inf h}\le\varepsilon\right\},
\]
and only count outer iterations, where $\mathbb K^\varepsilon_{\text{GNConvexArmijo}}(h)=0$ if $h(x_0)=\inf h$. 

The following result follows from the standard Armijo descent argument \citep{nocedal_numerical_2006}; its proof is included in Appendix \ref{app:proof_GNConvexArmijo}.

\begin{theorem}[Worst-Case Complexity of \eqref{alg:GNConvexArmijo}]\label{thm:WC_GNConvexArmijo}
    Let $h\in\mathcal F_{\mu,\underline\gamma}^{L,\bar\gamma}$ with $0<\mu\le L<+\infty$ and $0<\underline\gamma\le\bar\gamma<+\infty$, and let $(x_k)$ be the iterates of \eqref{alg:GNConvexArmijo} with $\beta,\tau\in(0,1)$. Then, for $\bar\alpha=\tau\min\{1,2(1-\beta)\underline\gamma^2/L\}$,
    \[
        \frac{h(x_k)-\inf h}{h(x_0)-\inf h}
        \le\left(1-\frac{2\mu\beta\bar\alpha}{\bar\gamma^2}\right)^k.
    \]
    Thus the relative accuracy complexity is bounded by, for $\varepsilon\in (0, 1)$,
    \[
        \sup_{h\in\mathcal F_{\mu,\underline\gamma}^{L,\bar\gamma}}\K^\varepsilon_{\text{GNConvexArmijo}}(h)
        \le\left\lceil\frac{\bar\gamma^2}{2\mu\beta\bar\alpha}\log\frac1\varepsilon\right\rceil.
    \]
\end{theorem}

The line search does not require knowledge of the values of $L$ and $\underline\gamma$, and instead adapts to them, thus preserving linear convergence. For fixed Armijo parameters, the smoothed complexity order of the fixed-step method is preserved.

\begin{theorem}[Smoothed Complexity of \eqref{alg:GNConvexArmijo}]\label{thm:SA_GNConvexArmijo}
    Let $0<\sigma,L,\bar\gamma<+\infty$ and $\varepsilon,\beta,\tau\in(0,1)$. We denote by $f_W$ a perturbation of $f$, parametrized by $\sigma$, as defined in Equation \eqref{eq:to_be_computed}. Then, as $n\to\infty$, the smoothed complexity of Algorithm \eqref{alg:GNConvexArmijo} over the function class $\mathcal F_{0,0}^{L,\bar\gamma}$ is bounded by
    \[
        \sup_{f\in\mathcal F_{0,0}^{L,\bar\gamma}}\ExpD{W}{\K^\varepsilon_{\text{GNConvexArmijo}}(f_W)}
        \le1+\mathcal O_{d,\beta,\tau}\left((1+\sigma^{-2}L)(1+\sigma^{-2}\bar\gamma^2)\log\frac1\varepsilon\right).
    \]
\end{theorem}
\begin{proof}
    Fix any $f\in\mathcal F_{0,0}^{L,\bar\gamma}$. By Theorem \ref{thm:WC_GNConvexArmijo},
    \[
        \K^\varepsilon_{\text{GNConvexArmijo}}(f_W)
        \le1+\frac1{2\beta\tau}\left(\frac{\bar\gamma_\sigma^2}{\mu_\sigma}
        +\frac{L_\sigma\bar\gamma_\sigma^2}{2(1-\beta)\mu_\sigma\underline\gamma_\sigma^2}\right)\log\frac1\varepsilon.
    \]
    The expectation of the first term is bounded by $\mathcal O_d(1+\sigma^{-2}\bar\gamma^2)$, as in the proof of Theorem \ref{thm:SA_GD1L}, by 
    \[
        \Exp{\frac{\bar\gamma_\sigma^2}{\mu_\sigma}}\le \Exp{\frac{\bar\gamma^2+\frac{\sigma^2s_{\max}(G)^2}{m}}{\frac{\sigma^2s_{\min}(G)^2}{m}}}=\Exp{\frac{m\sigma^{-2}\bar\gamma^2+s_{\max}(G)^2}{s_{\min}(G)^2}}\le \mathcal O_d(1+\sigma^{-2}\bar\gamma^2).
    \]
    The expectation of the second term is bounded as in the proof of Theorem \ref{thm:GNConvex_SA}. \qedhere
\end{proof}

\textbf{Discussion.} As in Section \ref{sec:GNConvex}, the smoothed complexity bound depends logarithmically on $\varepsilon$, thus providing a linear convergence rate rather than a sublinear rate as expected for merely convex functions. The dependency in the condition number remains square, as in Section \ref{sec:GNConvex}.

\begin{figure}[H]
    \centering
    \def\SAMethod{GNConvexArmijo}
    \input{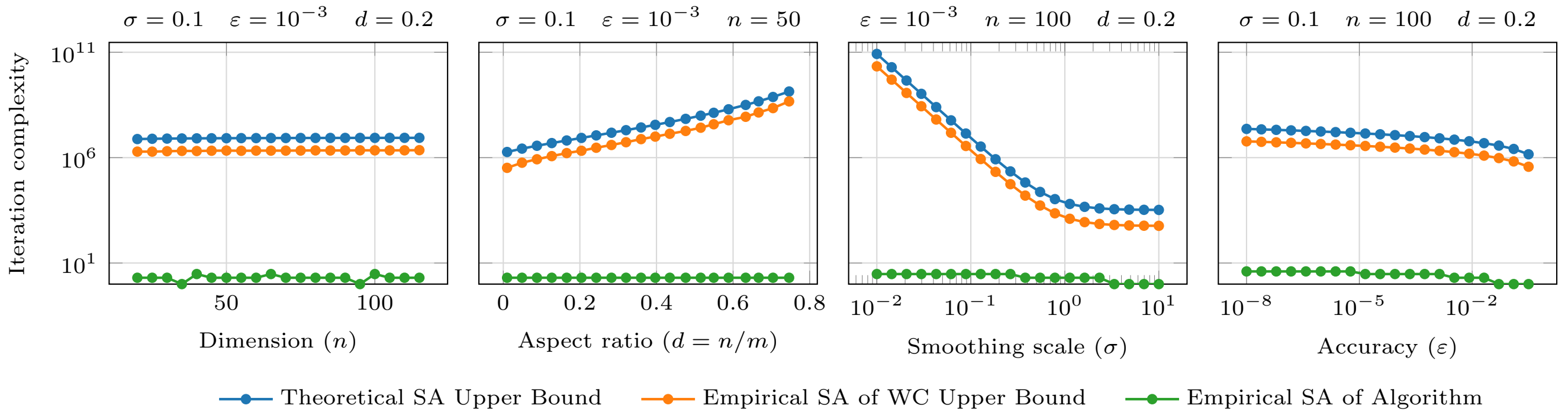}
    \caption{Smoothed complexity of Algorithm \eqref{alg:GNConvexArmijo} under Wishart perturbations.}
    \label{fig:GNConvexArmijo}
\end{figure}

\textbf{Numerical Illustration.} We use the same rank-one nonlinear least-squares objective as in Section \ref{sec:GNConvex}, namely
\[
    f(x)=\frac12\left(t(x)^2+\frac14\sin^2(t(x))\right),
    \qquad t(x)=\frac1{\sqrt n}\sum_{i=1}^n x_i,
\]
with the Armijo parameters $\beta=\tau=1/2$. The results are shown in Figure \ref{fig:GNConvexArmijo}. We note that our theory is overly pessimistic in this setting, where the empirical iteration count is always small, and our theory predicts a large value. This arises from the pessimistic worst-case bound in Theorem \ref{thm:WC_GNConvexArmijo}.

\subsection{Gauss-Newton Method - Nonconvex Case with Radial Monotonicity}\label{sec:GN_PL}

We now consider the nonconvex \textit{Nonlinear Least Squares} problem with Algorithm \eqref{alg:GN}, and further assume \textit{gradient domination}, as in Section \ref{sec:GD1L_PL}. Specifically, we denote by $\mathcal F_{\text{PL}, \eta, \underline{\gamma}}^{L, \bar \gamma}$ the subclass of $\mathcal F_{-\infty, \underline{\gamma}}^{L, \bar \gamma}$ that satisfies gradient domination with parameter $\eta>0$, where $\mathcal F_{-\infty, \underline{\gamma}}^{L, \bar \gamma}$ denotes the class of nonconvex $L$-smooth functions satisfying $\underline{\gamma}^2I\preceq J(x)^TJ(x)\preceq \bar\gamma^2 I$ for all $x$. We note that gradient domination implies that all local minimizers are global, and therefore we can measure complexity in terms of relative function value accurancy, namely, for $\varepsilon\in (0, 1)$,
\[
    \mathbb K^\varepsilon_{\text{GN, rad}}(h)=\min_{k\ge 0}\left\{k\colon \frac{h(x_k)-\inf h}{h(x_0)-\inf h}\le \varepsilon\right\},
\]
where $\mathbb K^\varepsilon_{\text{GN, rad}}(h)=0$ if $h(x_0)=\inf h$.

We first give the deterministic worst-case complexity under gradient domination, whose proof is postponed to Appendix \ref{app:proof_GNGradientDom}.

\begin{theorem}[Worst-Case Complexity of \eqref{alg:GN} under Gradient Domination]\label{thm:WC_GNGradientDom}
    Let $h\in \mathcal F_{\text{PL}, \eta, \underline{\gamma}}^{L, \bar \gamma}$ with $0<L, \eta< +\infty$ and $0< \underline{\gamma}\le \bar \gamma<+\infty$, and let $(x_k)$ be the iterates of \eqref{alg:GN}. Then,
    \[
        \frac{h(x_k)-\inf h}{h(x_0)-\inf h}\le \left(1-\frac{\eta \underline{\gamma}^2}{L\bar\gamma^2}\right)^k.
    \]
    Thus, the relative function value accuracy complexity is bounded by, for $\varepsilon\in (0, 1)$,
    \[
        \sup_{h\in \mathcal F_{\text{PL}, \eta, \underline{\gamma}}^{L, \bar \gamma}}\mathbb K^\varepsilon_{\text{GN, rad}}(h)\le \left\lceil\frac{L\bar\gamma^2}{\eta \underline{\gamma}^2}\log \frac{1}{\varepsilon}\right\rceil.
    \]
\end{theorem}
As in Sections \ref{sec:GD1L_PL} and \ref{sec:CNM_PL}, gradient domination may be covered through perturbation of radially monotone functions. We denote by $\mathcal F_{\text{rad}}^{L, \bar \gamma}(\tilde x)$ the subclass of $\mathcal F_{-\infty, 0}^{L, \bar \gamma}$ that satisfies radial monotonicity. We now state the smoothed complexity guarantee.
\begin{theorem}[Smoothed Complexity of \eqref{alg:GN} under Radial Monotonicity]
    Let $0<\sigma, L, \bar \gamma <+\infty$, $\varepsilon\in (0, 1)$, and let $\tilde x\in \R^n$. We denote by $f_W$ a perturbation of $f$, parametrized by $\sigma$, as defined in Equation \eqref{eq:to_be_computed}. Then, as $n\to \infty$, the smoothed complexity of Algorithm \eqref{alg:GN} over the function class $\mathcal F_{\text{rad}}^{L, \bar \gamma}(\tilde x)$ is bounded by 
    \[
        \sup_{f\in \mathcal F_{\text{rad}}^{L, \bar \gamma}(\tilde x)}\ExpD{W}{\mathbb K^\varepsilon_{\text{GN, rad}}(f_W)}\le 1+\mathcal O_d\left((1+\sigma^{-2}L)^2(1+\sigma^{-2}\bar\gamma^2)\log\frac{1}{\varepsilon}\right).
    \]
\end{theorem}
\begin{proof}
    Fix any $f\in \mathcal F_{\text{rad}}^{L, \bar \gamma}(\tilde x)$. Theorem \ref{thm:WC_GNGradientDom} and the preceding bounds give 
    \[
        \mathbb K^\varepsilon_{\text{GN, rad}}(f_W)\le1+ \frac{(L+\sigma^2\lambda_{\max}(W))^2(\bar\gamma^2+\sigma^2\lambda_{\max}(W))}{\sigma^6\lambda_{\min}(W)^3}\log\frac{1}{\varepsilon}.
    \]
    Expanding the numerator, applying Cauchy-Schwarz to each mixed moment and using Theorem \ref{thm:moments} proves the result.
\end{proof}

\textbf{Discussion.} We note that, akin to Gradient Descent in Section \ref{sec:GD1L_PL}, the radial monotonicity setting adds a factor $1+\sigma^{-2}L$ to the smoothed bound, compared to the convex case in Section \ref{sec:GNConvex}.

\textbf{Numerical Illustration.} We consider the following nonlinear least squares problem
\[
    f(x)=\frac{1}{2}\|r(x)\|^2=\frac{1}{2}\sum_{i=1}^{\lfloor n/2\rfloor}r_i(x)^2, \quad r_i(x)=\frac{x_i}{\sqrt{1+x_i^2}}\quad \text{for $i=1, \ldots, \lfloor n/2\rfloor$}, \quad \tilde x=0.
\]
We note that $(x-\tilde x)^T\nabla f(x)=\sum_i x_i^2/(1+x_i^2)^2\ge 0$ and $f(\tilde x)=\min f$, so that the function is radially monotone. We compute that
\[
    \partial_{ii}f(x)=\frac{1-3x_i^2}{(1+x_i^2)^3}\quad\text{and}\quad J_ii(x)=(1+x_i^2)^{-3/2},
\]
so that the objective is nonconvex, globally smooth with $L=1$ and has Jacobian singular values bounded by $\bar \gamma=1$. The Jacobian is rank-deficient at every $x$, as the sum only spans the first $\rfloor n/2\rfloor$ coordinates, and does not satisfy gradient domination as, for instance, along $x=te_1$, the objective tends to a finite constant as $t\to \infty$. 

\begin{figure}[H]
    \centering
    \def\SAMethod{GNSaturating}
    \input{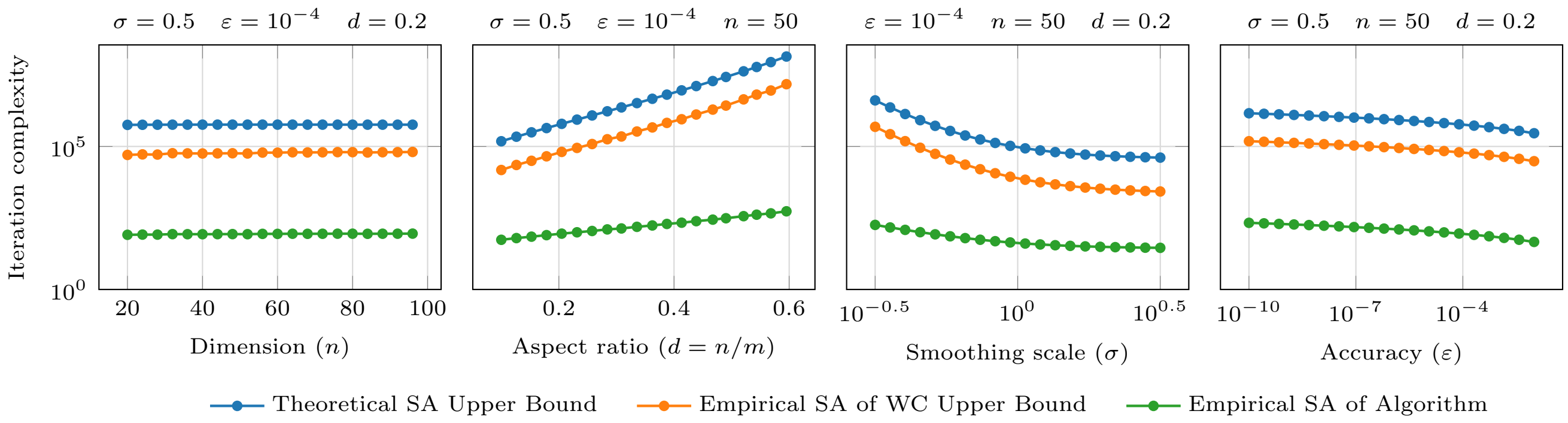}
    \caption{Smoothed complexity of Algorithm \eqref{alg:GN} under Wishart perturbations.}
    \label{fig:GN_PL}
\end{figure}

The results are shown in Figure \ref{fig:GN_PL}. As in Sections \ref{sec:GD1L_PL} and \ref{sec:CNM_PL}, we observe that our bounds are more pessimistic in its dependencies that in the convex counterparts, but that the trends are correct. 

\section{Conclusion}

In this paper, we developed a smoothed analysis framework for a range of classical algorithms in continuous optimization, showing that small random perturbations can lead to substantially sharper complexity guarantees than those suggested by a worst-case analysis. Specifically, we apply methods suited for non-degenerate functions to small perturbations of degenerate functions, and recover the faster best-known rates for the non-degenerate class. The resulting bounds obtained in the framework help explain why algorithms that may exhibit poor worst-case guarantees nevertheless can perform well on nearby instances.

\paragraph{Outlook.} Our results suggest that smoothed analysis can provide a useful intermediate perspective between worst-case and average-case analysis for continuous optimization methods. Several questions however remain open. Firstly, the current analysis heavily relies on worst-case results, which inevitably weakens the smoothed complexity results obtained. Secondly, it would be interesting to extend the framework to more nonconvex settings, beyond the radial monotonicity setting in Sections \ref{sec:GD1L_PL}, \ref{sec:CNM_PL} and \ref{sec:GN_PL}. Finally, all algorithms in the present work are deterministic. An open question remains how to tackle stochastic algorithms, and specifically the algorithm of \cite[Section 3.3]{cortild_stochastic_2026}, which was one motivation of this work.

\paragraph{Acknowledgments.} Daniel Cortild acknowledges the support of the Clarendon Funds Scholarships. Coralia Cartis' work was supported by the Hong Kong Innovation and Technology Commission (InnoHK Project CIMDA) and by the EPSRC grant EP/Y028872/1, Mathematical Foundations of Intelligence: An “Erlangen Programme” for AI.

{
\small
\setlength{\bibsep}{0pt}
\bibliographystyle{apalike}
\bibliography{references}
}

\appendix

\section{Proofs}

\subsection{Proof of Theorem \ref{thm:GNConvex}}\label{app:proof_GN_Convex}


Let $B_k=J_k^TJ_k$. Since $\underline\gamma^2I\preceq B_k\preceq\bar\gamma^2I$,
\[
    g_k^Ts_k=-g_k^TB_k^{-1}g_k\le-\bar\gamma^{-2}\|g_k\|^2,
    \qquad
    \|s_k\|^2=g_k^TB_k^{-2}g_k\le-\underline\gamma^{-2}g_k^Ts_k.
\]
By $L$-smoothness, for any $\alpha>0$,
\begin{equation}\label{eq:GN_decrease}
    h(x_k+\alpha s_k)
    \le h(x_k)+\alpha g_k^Ts_k+\frac{L\alpha^2}{2}\|s_k\|^2 \le h(x_k)+\alpha\left(1-\frac{L\alpha}{2\underline\gamma^2}\right)g_k^Ts_k.
\end{equation}
The guaranteed decrease is maximized by the fixed step size $\alpha=\underline\gamma^2/L$. Hence,
\[
    h(x_{k+1})-\min h
    \le h(x_k)-\min h-\frac{\underline\gamma^2}{2L\bar\gamma^2}\|g_k\|^2
    \le \left(1-\frac{\mu\underline\gamma^2}{L\bar\gamma^2}\right)(h(x_k)-\min h),
\]
where the last inequality uses $\|g_k\|^2\ge2\mu(h(x_k)-\min h)$. Iterating proves the claim.

\subsection{Proof of Theorem \ref{thm:WC_GNConvexArmijo}}\label{app:proof_GNConvexArmijo}

Let $B_k=J_k^TJ_k$. As in Appendix \ref{app:proof_GN_Convex},
\[
    g_k^Ts_k\le-\bar\gamma^{-2}\|g_k\|^2,
    \qquad \|s_k\|^2\le-\underline\gamma^{-2}g_k^Ts_k.
\]
By Equation \eqref{eq:GN_decrease}, Armijo's condition is satisfied whenever $\alpha\le{2(1-\beta)\underline\gamma^2}/{L}$. Therefore, $\alpha_k\ge\bar\alpha=\tau\min\{1,2(1-\beta)\underline\gamma^2/L\}$. By strong convexity,
\begin{align*}
    h(x_{k+1})-\min h
    \le h(x_k)-\min h+\beta\alpha_kg_k^Ts_k 
    \le h(x_k)-\min h-\beta\bar\alpha\bar\gamma^{-2}\|g_k\|^2 
    \le\left(1-\frac{2\mu\beta\bar\alpha}{\bar\gamma^2}\right)(h(x_k)-\min h).
\end{align*}
Iterating and using $1-t\le e^{-t}$ proves the two claims.



\subsection{Proof of Theorem \ref{thm:WC_GNGradientDom}}\label{app:proof_GNGradientDom}

Taking $\alpha=\underline\gamma^2/L$ in Equation \eqref{eq:GN_decrease} gives
\begin{equation*}\label{eq:something_important}
    h(x_{k+1})\le h(x_k)-\frac{\underline\gamma^2}{2L\bar\gamma^2}\|g_k\|^2.
\end{equation*}
By gradient domination, we obtain
\[
    h(x_{k+1})-\inf h\le h(x_k)-\inf h-\frac{\underline\gamma^2}{2L\bar\gamma^2}\|g_k\|^2\le h(x_k)-\inf h-\frac{\eta\underline\gamma^2}{L\bar\gamma^2}(h(x_k)-\inf h),
\]
which proves the first inequality.
\section{Explicit Constants}\label{app:constants}

We sharpen both moment constants in Theorem \ref{thm:moments}. This does not affect any of the results in the main body, which only require bounds uniform in $m$ and $n$.

\begin{proposition}[Positive Moments of the Largest Singular Value]\label{prop:moment_max_explicit}
    Let $q>0$. If $q\ge 1$, then
    \[
        \Exp{m^{-q/2}s_{\max}(G)^q}
        \le \left(1+\sqrt{\frac nm}+\frac{c_q}{\sqrt m}\right)^q,
        \qquad
        c_q=\left(q2^{q/2}\Gamma(q/2)\right)^{1/q}.
    \]
    If $0<q<1$, the same conclusion holds with $c_q=\sqrt{2\pi}$. In particular, whenever $n/m\le d<1$,
    \[
        C^{\max}_{n,m,q}=\left(1+\sqrt{\frac nm}+\frac{c_q}{\sqrt m}\right)^q
        \le C^{\max}_{d,q}=\left(1+\sqrt d+c_q\right)^q.
    \]
\end{proposition}
\begin{proof}
    Let
    \[
        Z_+=\left[s_{\max}(G)-\sqrt m-\sqrt n\right]_+.
    \]
    The Gaussian tail bound \cite[Equation (2.3)]{rudelson_nonasymptotic_2010} gives $\P(Z_+\ge t)\le 2e^{-t^2/2}$ for every $t\ge0$. Hence, for any $q>0$,
    \begin{align*}
        \Exp{Z_+^q}
        &=q\int_0^\infty t^{q-1}\P(Z_+\ge t)dt \le 2q\int_0^\infty t^{q-1}e^{-t^2/2}dt
        =q2^{q/2}\Gamma(q/2).
    \end{align*}
    For $q\ge1$, Minkowski's inequality therefore yields
    \[
        \left\|m^{-1/2}s_{\max}(G)\right\|_q
        \le 1+\sqrt{\frac nm}+\frac{\|Z_+\|_q}{\sqrt m}
        \le 1+\sqrt{\frac nm}+\frac{c_q}{\sqrt m}.
    \]
    For $0<q<1$, Jensen's inequality reduces the claim to the case $q=1$, for which $c_1=\sqrt{2\pi}$. \qedhere
\end{proof}

\begin{proposition}[Negative Moments of the Smallest Singular Value]\label{prop:moment_min_explicit}
    Let $p=(m-n+1)/2$ and $0<q<2p$, and define
    \[
        a_{n,m}=\frac{n2^{-p}\Gamma(3/2)\Gamma(m/2+1/2)}
        {p\Gamma(n/2+1)\Gamma(p)\Gamma(p+1/2)}.
    \]
    Then
    \[
        \Exp{m^{q/2}s_{\min}(G)^{-q}}
        \le C^{\min}_{n,m,q}=
        \frac{2p}{2p-q}\left(ma_{n,m}^{1/p}\right)^{q/2}.
    \]
    Moreover, if $n/m\le d<1$ and $m,n\to\infty$ proportionally, then, for all sufficiently large $n$,
    \[
        C^{\min}_{n,m,q}\le C^{\min}_{d,q}= 2K_d^q,
        \qquad
        K_d=\left(\frac{ed^{1-1/(1-d)}}{(1-d)^2}\right)^{1/2}.
    \]
\end{proposition}
\begin{proof}
    Let $0\le\lambda_1\le\cdots\le\lambda_n$ be the ordered eigenvalues of $GG^T$. From their joint density \citep{dumitriu_matrix_2002}, for $c_{n, m}$ the ordered-density constant,
    \begin{align*}
        \P(\lambda_1\le s)
        &\le \frac{c_{n,m}}{c_{n-1,m+1}}\int_0^s\lambda^{p-1}d\lambda
        =a_{n,m}s^p,
    \end{align*}
    where the displayed expression for $a_{n,m}$ follows by simplifying the ratio of normalizing constants. Hence,
    \begin{align*}
        \Exp{\lambda_1^{-q/2}}
        &=\frac q2\int_0^\infty t^{q/2-1}\P(\lambda_1\le t^{-1})dt 
        \le \frac q2\int_0^{a_{n,m}^{1/p}}t^{q/2-1}dt
        +\frac q2a_{n,m}\int_{a_{n,m}^{1/p}}^\infty t^{q/2-p-1}dt 
        =\frac{2p}{2p-q}a_{n,m}^{q/(2p)}.
    \end{align*}
    Since $s_{\min}(G)^2=\lambda_1$, multiplying by $m^{q/2}$ proves the first claim. Finally, the Stirling bounds for the Gamma function give
    \[
        a_{n,m}\le \left(\frac{K_d^2}{m}\right)^p n^{1/2}p^{-1} \implies 
        C^{\min}_{n,m,q} \le \frac{2p}{2p-q}K_d^q,
    \]
    where the implications holds once $p\ge \sqrt n$, and the stated uniform constant follows once $p\ge q$. \qedhere
\end{proof}

\paragraph{From Theorem \ref{thm:moments}.} Propositions \ref{prop:moment_max_explicit} and \ref{prop:moment_min_explicit} give the constants
\[
    C^{\max}_{n,m,q}=\left(1+\sqrt{\frac nm}+\frac{c_q}{\sqrt m}\right)^q\le C^{\max}_{d,q}=(1+\sqrt d+c_q)^q,
    \qquad
    c_q=\begin{cases}
        \left(q2^{q/2}\Gamma(q/2)\right)^{1/q}, & q\ge1,\\
        \sqrt{2\pi}, & 0<q<1,
    \end{cases}
\]
\[
    C^{\min}_{n,m,q}=\frac{2p}{2p-q}\left(ma_{n,m}^{1/p}\right)^{q/2}\le 
    C^{\min}_{d,q}=2\left(\frac{ed^{1-1/(1-d)}}{(1-d)^2}\right)^{q/2}.
\]

\paragraph{From Smoothed Complexity Theorems.} 

For $u,v>0$, $t\geq0$, and $r>1$, define
\[
    p=\frac{m-n+1}{2}\qquad  \frac nm\le d<1, \qquad M_{u,v}=\sqrt{C^{\max}_{d,4u}C^{\min}_{d,4v}},
    \qquad\text{and}\qquad
    \Phi_r(t)=\log_r\left(1+\log_2(1+t)\right).
\]
Define the auxiliary constants
{\small\[
\begin{cases}
    C_{\text{GD}}
    &=\sigma^{-2}LC^{\min}_{d,2}+M_{1,1},\\
    C_{\text{GD,rad}}
    &=\sigma^{-4}L^2C^{\min}_{d,4}+2\sigma^{-2}LM_{1,2}+M_{2,2},\\
    C_{\text{acc}}
    &=\sigma^{-1}\sqrt L C^{\min}_{d,1}+M_{1/2,1/2},\\
    C_{\text{CNM}}
    &=6.25\cdot 2^{3/4}L_H^{1/2}\left(
       \Delta^{1/4}\sigma^{-3/2}C^{\min}_{d,3/2}
       +2^{-1/4}D^{1/2}\sigma^{-1}M_{1/4,3/4}
       \right),\\
    C_{\text{CNM,rad}}
    &=8\cdot 2592^{1/4}L_H^{1/2}\Bigl[
       \Delta^{1/4}\sigma^{-3/2}
       \left(\sigma^{-3/2}L_G^{3/4}C^{\min}_{d,3}+M_{3/4,3/2}\right)+2^{-1/4}D^{1/2}\sigma^{-1}
       \left(\sigma^{-3/2}L_G^{3/4}M_{1/4,3/2}+M_{1,3/2}\right)
       \Bigr],\\
    C_{\text{Hyb, A}}
    &=2L_H^2\left(
       4L_G^2\Delta^2\sigma^{-16}C^{\min}_{d,16}
       +(4\Delta^2+L_G^2D^4)\sigma^{-12}M_{2,8}
       +D^4\sigma^{-8}M_{4,8}
       \right)^{1/2},\\
    C_{\text{Hyb}}
    &=\left(2L_G^2\sigma^{-4}C^{\min}_{d,4}
       +2M_{2,2}\right)^{1/2}\bigl(\log(1+C_{\text{Hyb, A}})+1\bigr),\\
    C_{\text{GN}}
    &=\sigma^{-4}L \bar\gamma^2C^{\min}_{d,4}+\sigma^{-2}(L+\bar\gamma^2)M_{1,2}+M_{2,2},\\
    C_{\text{GN,Armijo}}
    &=\frac1{2\beta\tau}\left(
       \sigma^{-2}\bar\gamma^2C^{\min}_{d,2}+M_{1,1}
       +\frac{C_{\text{GN}}}{2(1-\beta)}\right),\\
    C_{\text{GN,rad}}
    &=\sigma^{-4}L^2\sigma^{-2}\bar\gamma^2C^{\min}_{d,6}
       +\sigma^{-4}(L^2+2L\bar\gamma^2)M_{1,3}
       +\sigma^{-2}(2L+\bar\gamma^2)M_{2,3}+M_{3,3}. \\
\end{cases}
\]}The explicit bounds are given by
\[
\begin{cases}
    \Exp{\K^\varepsilon_{\text{GD}}(f_W)}
    &\le1+C_{\text{GD}}\log\frac1\varepsilon,
       \quad\text{for }p\geq\max\{\sqrt n,4\},\\
    \Exp{\K^\varepsilon_{\text{GD,rad}}(f_W)}
    &\le1+C_{\text{GD,rad}}\log\frac1\varepsilon,
       \quad\text{for }p\geq\max\{\sqrt n,8\},\\
    \Exp{\K^\varepsilon_{\text{NAG}}(f_W)}
    &\le1+C_{\text{acc}}\log\frac2\varepsilon,
       \quad\text{for }p\geq\max\{\sqrt n,2\},\\
    \Exp{\K^\varepsilon_{\text{CG}}(f_W)}
    &\le
       1+\frac{C_{\text{acc}}}{2}\log\frac4\varepsilon,
       \quad\text{for }p\geq\max\{\sqrt n,2\},\\
    \Exp{\K^\varepsilon_{\text{CNM}}(f_W)}
    &\le2+C_{\text{CNM}}
       +\Phi_2\left({C^{\max}_{d,6}\sigma^6}/(2L_H^2\varepsilon)\right),
       \quad\text{for }p\geq\max\{\sqrt n,3\},\\
    \Exp{\K^\varepsilon_{\text{CNM,rad}}(f_W)}
    &\le2+C_{\text{CNM,rad}}
       +\Phi_{4/3}({C^{\max}_{d,6}\sigma^6}/(2592L_H^2\varepsilon)),
       \quad\text{for }p\geq\max\{\sqrt n,6\},\\
    \Exp{\K^\varepsilon_{\text{Hyb}}(f_W)}
    &\le3+C_{\text{Hyb}}
       +\Phi_2\left({C^{\max}_{d,6}\sigma^6}/(2L_H^2\varepsilon)\right),
       \quad\text{for }p\geq\max\{\sqrt n,32\},\\
    \Exp{\K^\varepsilon_{\text{GNConvex}}(f_W)}
    &\le1+C_{\text{GN}}\log\frac1\varepsilon,
       \quad\text{for }p\geq\max\{\sqrt n,8\},\\
    \Exp{\K^\varepsilon_{\text{GNConvexArmijo}}(f_W)}
    &\le1+C_{\text{GN,Armijo}}\log\frac1\varepsilon,
       \quad\text{for }p\geq\max\{\sqrt n,8\},\\
    \Exp{\K^\varepsilon_{\text{GN, rad}}(f_W)}
    &\le1+C_{\text{GN,rad}}\log\frac1\varepsilon,
       \quad\text{for }p\geq\max\{\sqrt n,12\}.
\end{cases}
\]

\end{document}